\documentclass[10pt,psamsfonts]{amsart}
\usepackage{a4}
\usepackage{tikz}
\usetikzlibrary{matrix}
\usepackage{indentfirst}
\usepackage{graphics,amsthm}
\usepackage{color}
\usepackage{epsfig}
\usepackage{fancybox}
\usepackage{mathrsfs}
\usepackage{hyperref}
\usepackage{newlfont}
\usepackage{amsmath,amsfonts,amstext,amssymb,bezier,amsthm}
\usepackage{subfigure}
\usepackage[active]{srcltx}
\usepackage[all]{xy}
\usepackage[english]{babel}
\usepackage[latin1]{inputenc}
\usepackage{graphicx,psfrag}
\usepackage{fancyhdr}
\usepackage{niceframe}
\usepackage{xspace,setspace} %sum\'{a}rio
\usepackage{dsfont}
\usepackage{caption}%[font=small,labelfont=bf]{caption}

\newtheorem{theorem}{Theorem}%[section]
\newtheorem{proposition}{Proposition}%[section]
\newtheorem{corollary}{Corollary}%[section]
\newtheorem{lemma}{Lemma}%[section]
\newtheorem{definition}{Definition}%[section]
\theoremstyle{definition}
\newtheorem{example}{Example}%[section]
\theoremstyle{remark}
\newcommand{\prooff}{\noindent \textbf{Proof: }}
\newcommand{\cqd}{{{\hfill $\square$}}\\}

\newcommand{\Z}{\mathbb{Z}}

\newcommand{\I}{\textbf{I}}
\newcommand{\J}{\textbf{J}}
\newcommand{\K}{\textbf{K}}
\newcommand{\PP}{\textbf{P}}

\begin{document}

\renewcommand{\arraystretch}{1.5}

\author[Franzosa]
{R. Franzosa$^1$}

\author[de Rezende]{K. A. de Rezende$^2$}
\thanks{$^2$Partially supported by CNPq under grant 302592/2010-5 and FAPESP under grant 2012/18780-0}

\author[Vieira]{E. R. Vieira$^3$}
\thanks{$^3$Supported by FAPESP under grant 2010/19230-8.}

\address{$^1$ Departament of Mathematics and Statistics,
University of Maine, Orono, Maine, USA} \email{robert$\_$franzosa@umit.maine.edu}

\address{$^2$ Departamento de Matematica, Universidade
Estadual de Campinas, 13083--859, Campinas, SP,
Brazil} \email{ketty@ime.unicamp.br}

\address{$^3$ Departamento de Matematica, Universidade
Estadual de Campinas, 13083--859, Campinas, SP,
Brazil} \email{ewertonrvieira@gmail.com}

\title[Generalized Topological Transition Matrix.]
{Generalized Topological Transition Matrix.}

\subjclass[2010]{Primary 37B30; 37D15;  Secondary 70K70; 70K50; 55T05}

\keywords{Conley index, Connection Matrices, Transition Matrices, Morse-Smale System; Sweeping Method; Spectral Sequence}

\maketitle

\begin{abstract}
This article represents a major step in the unification of the theory of algebraic, topological and singular transition matrices by introducing a definition which is a generalization that encompasses all of the previous three. 
When this more general transition matrix satisfies  the additional requirement that it covers flow-defined Conley-index isomorphisms, one proves algebraic and connection-existence properties. These general transition matrices with this covering property are referred to as generalized topological transition matrices and are used to consider connecting orbits of Morse-Smale flows without periodic orbits, as well as those in a continuation associated to a dynamical spectral sequence.

%In this article this generalized definition This is accomplished, initially by presenting a generalization of the topological transition matrix. The main defining condition for these new matrices is that they cover flow-defined Conley-index isomorphisms. In contrast to the classical case, there is no need to require that there are no connections at the initial and final parameters of a continuation. 
%Algebraic and connection-existence properties for the generalized topological transition matrix that extend those of the classical case are proved herein. These results are used to consider connecting orbits of Morse-Smale flows without periodic orbits, as well as those in a continuation associated to a dynamical spectral sequence.

%In this article we present a generalization of the topological transition matrix, by removing the requirement of no connections at the initial and final parameters of a continuation. We prove similar properties for generalized topological transition matrix that extend the ones in the classical case. We also apply these results to Morse-Smale flows without periodic orbits, as well as to a continuation associated to a dynamical spectral sequence.
\end{abstract}

%$$
%\mathscr{ABCDEFGHIJLMNOPQRSTUVXZWY}
%$$
%$$
%\mathfrak{ABCDEFGHIJLMNOPQRSTUVXZWY}
%$$
%$$
%\mathcal{ABCDEFGHIJLMNOPQRSTUVXZWY}
%$$

%\tableofcontents

\section{Introduction}

A challenging question in the study of dynamical
 systems is that of the existence of global bifurcations. The difficulty in detecting such bifurcation orbits is the
fact that one must analyze the dynamical system globally. Topological techniques for global analysis are, therefore, a perfect fit for such an investigation.
In particular,  Conley index theory has proven to be quite useful in this role, as can be seen by the ample use of connection and transition matrices in bifurcation-related results. See \cite{C}, \cite{CF}, \cite{FiM}, \cite{Fr1}, \cite{Fr2}, \cite{Fr3} and \cite{MM1}.

%A challenging question in the study of dynamical
%systems is that of the existence of global bifurcations. The difficulty in detecting such bifurcation orbits is the
%fact that one must analyze the dynamical system globally thus making topological techniques a perfect fit.
%In particular, the usefulness of Conley index theory in this role has been demonstrated by the ample use of connection matrices in bifurcation related results. See \cite{C}, \cite{Fr1}, \cite{Fr2} and \cite{Fr3}.

Connection matrices have been extensively studied and can be computed by numerical techniques \cite{BM}, \cite{BR} and \cite{E}. Their continuation properties have proven useful in detecting global bifurcations. In particular, the continuation theorem \cite{Fr3}
states that the connection matrices
of an admissible ordering are invariant under local continuation. Yet, under global continuation, sets of connection matrices can undergo change. For instance, if there is a continuation between parameters with unique but different connection matrices, then within the continuation there must be a parameter value with nonunique connection matrices. At such a parameter value the system typically has a global bifurcation.

%These matrices have been extensively studied and can be computed by numerical techniques \cite{BM-H}. However, the lack of uniqueness of  connection matrices of certain dynamical systems seemed to indicate the presence of global bifurcation. The continuation theorem \cite{Fr3}
%states that the connection matrices
%of an admissible ordering are invariant under local continuation. Thus, at parameter
%values where bifurcations occur, one would expect nonuniqueness of connection
%matrices.

In other words, Morse decompositions and connection matrices provide a supporting structure within which  global bifurcations can be detected, particularly via changes in the associated algebraic structures. These differences that occur in connection matrices under continuation, which can naturally be identified algebraically, 
was the main motivation for the introduction of transition matrices as a combinatorial mechanism to keep track of these changes.
%and tracking these changes was the main motivation for the introduction of transition matrices. 
These transition matrices have since appeared in the literature  under several guises: singular \cite{R}, topological \cite{MM1}, and algebraic \cite{FM}. These three types of matrices are defined differently (particularly under contrasting conditions) and have distinct properties. On the other hand, due to underlying similarities in the definitions and their corresponding properties,  a unified theory for transition matrices has long been called for.

%Morse decompositions and connection matrices provide a supporting structure within which one can detect global bifurcations. The appropriate framawork within Conley theory to detect these bifurcations was the main motivation for the introduction of transition matrices. These matrices appear under several guises, singular \cite{R}, topological \cite{MM1} and algebraic \cite{FM}. These three different constructions of transition matrices, that appear in the literature are defined differently and have different properties.

In this paper we briefly introduce the generalization which unifies the theory. We focus on an initial and important step toward understanding the properties of this newly defined and more general  transition matrix,  which has the additional property  that it covers flow-defined Conley-index isomorphisms. We refer to these matrices as generalized topological transition matrices and prove several properties they possess. In contrast to the classical case, we do not require that there are no connections at the initial and final parameters of a continuation.

%For instance, in \cite{FM}, is is proved that algebraic transition matrices can be used to generate all connection matrices at a particular parameter value. The consequence of this is that the burden to find changes in global dynamics is shifted to transition matrices.

%By undertaking the task of unifying the theory of transition matrix we actually present a framework which makes possible the defining of a generalized transition matrix which encompasses all transition matrices previously defined, such as those in \cite{R}, \cite{MM1} and \cite{FM}.

%However, in order to prepare the foundation needed to unify the transition theory, in this paper we generalize the classical topological transition matrix and explore its properties. To achieve this one must get rid of the restriction of no connections between the Morse sets at the end parameter values in a continuation.

In section 2, we establish properties of the generalized topological transition matrices - including connecting orbit existence results - corresponding to those of the classical topological transition matrix. In section 3, we apply this new theory to Morse-Smale flows without periodic orbits. In this setting one demonstrates uniqueness and provides a simple way to compute the generalized topological transition matrix. In the last section, we see how the generalized topological transition matrices can be obtained from a continuation associated to a dynamical spectral sequence.

%In section 2, we establish properties of the generalized topological transition matrix - including connecting orbit existence results - corresponding to those of the classical topological transition matrix. Moreover, in section 3, we apply this new theory to Morse-Smale flows without periodic orbits. In this setting, one establishes uniqueness and an easier way to compute the generalized topological transition matrix. In the last section, we see how generalized topological transition matrices can be obtained from a continuation associated to a dynamical spectral sequence.

We assume that the reader is familiar with the basic ideas in Conley Index Theory, including Morse decompositions, homology index braids, connection matrices, etc. (see \cite{C}, \cite{Fr1}, \cite{Fr2}, \cite{Fr3}, \cite{MMr} and \cite{S}).

%In classical Morse theory the gradient vector field of a Morse function on a compact manifold gives rise to a flow with a finite number of hyperbolic fixed points. C. Conley introduced a generalization of this in the form of a Morse decomposition. Morse theory comes from its ability to relate local information (the Morse indices of the hyperbolic critical points) to global information
%(the homology of the compact manifold). This relation appears in the form of the
%Morse inequalities. In the Conley theory, the Morse inequalities are replaced by connection
%matrix.

Let $\varphi$ be a continuous flow on a locally compact Hausdorff space and let $S$ be a compact invariant set under $\varphi$. A \textit{Morse decomposition} of $S$ is a collection of mutually disjoint compact invariant subsets of $S$,
$$
\mathcal{M}(S)=\{M(p)\ |\ \pi\in\PP\}
$$
indexed by a finite set $\PP$, where each set $M(p)$ is called a \textit{Morse set}. A partial order $<$ on $\PP$ is called \textit{admissible ordering} if for $x\in S\backslash \bigcup_{\pi\in\PP}M(p)$ there exists $p<q$ such that $\alpha(x)\subseteq M(p)$ and $\omega(x)\subseteq M(q).$ The flow defines an
admissible ordering of $M$, called the \textit{flow ordering} of $M$, denoted $<_F$, and such
that $M(\pi)<_F M(\pi') $ if and only if there exists a sequence of distinct elements
of $P: \pi = \pi_0,\ldots,\pi_n=\pi'$, where $C(M(\pi_j),M(\pi_{j-1}))$, the set of connecting orbit between $M(\pi_j)$ and $M(\pi_{j-1})$, is nonempty for each $j = 1, \ldots ,n$. Note that every admissible ordering of $M$ is an extension of $<_F$.

In the Conley theory one begins with the Conley index for {isolated invariant sets}, i. e.,  $S\subseteq X$ is an  \textit{isolated invariant set} if there exists a compact set $N\subseteq X$ such that
$$
S=Inv(N,\varphi)=\{x\in N|\ \mathcal{O}(x) \subseteq N\}.
$$
The \textit{homological Conley index of} $S$, $H_\ast(S)$ is the homology of the pointed space $(N\backslash L)$, where $(N,L)$ is an index pair for $S$. Setting $$M(\I)=\bigcup_{\pi\in\I}M(\pi)\cup\bigcup_{\pi,\pi'\in \I}C(M(\pi'),M(\pi)),$$ the Conley index of $M(\I)$, $CH_\ast(M(\I))$, in short $H_\ast(\I)$, is well defined, since $M(\I)$ is an isolated invariant set for all $\I\in\I(<)$.

Given $\mathcal{M}(S)$, a Morse decomposition of $S$, the existence of an admissible ordering on $\mathcal{M}(S)$ implies that any recurrent dynamics in $S$ must be contained within the Morse sets, thus the dynamics off the Morse sets must be gradient-like. For this reason, Conley index theory refers to the dynamics within a Morse set as local dynamics and off the Morse sets as global dynamics. We briefly introduce the connection matrix theory, which addresses this latter aspect.

\begin{definition}\label{connection_matrix}
Given $\mathcal{G}$, a graded module braid over $<$, and $\mathcal{C} = \{C(\pi)\}_{\pi\in\PP}$,
a collection of graded modules, let $\Delta: \bigoplus_{\pi\in\PP} C(\pi) \rightarrow \bigoplus_{\pi\in\PP} C(\pi)$ be a $<$-upper
triangular boundary map. If $ \mathcal{H}\Delta$, the graded module braid generated by $\Delta$, is isomorphic to $\mathcal{G}$, and $C(p)$ is isomorphic to $\mathcal{G}(p)$ then $\Delta$
is called a connection matrix of $\mathcal{G}$.
\end{definition}

To simplify notation, for $\I \in \I(<)$ we denote $\bigoplus_{\pi\in\I}C(\pi)$ by $C(\I)$, and the
corresponding homology module in $\mathcal{H}\Delta$ by $H(\I)$. In particular, the homology index braid of an admissible ordering of a Morse decomposition $\mathcal{G}=\{H_\ast(\I)\}_{\I\in\I(<)}$ is an example of a graded module braid. In this setting  a $<$-upper
triangular boundary map
$$\Delta:\displaystyle \bigoplus_{\pi\in \PP}CH_\ast (M(\pi)) \rightarrow \displaystyle \bigoplus_{ \pi\in  \PP}CH_{\ast-1} (M( \pi))$$
satisfying Definition \ref{connection_matrix} for $\mathcal{C}\Delta=\{CH_\ast(M(\pi))\}_{\pi\in\PP}$ is called the \textit{connection matrix for  a Morse decomposition}. Since in this paper, our aim is to work with topological transition matrices, we focus on connection matrices for Morse decompositions. Thus, let $\mathcal{CM}(<)$ denote the set of all connection matrices for a given ($<$-ordered) Morse decomposition $\mathcal{M}(S)$.

{One of the key features in Conley theory is its invariance under continuation. Since the connection matrices for Morse decompositions
are algebraically derived from the homology Conley index braid, this seems to indicate that connecting orbits that persist over open sets in parameter space are identified by connection matrices. We now define Conley index continuation.}

Let $\Gamma$ be a Hausdorff topological space, $\Lambda$ a compact, locally contractible, connected metric space and $X$ a locally compact metric space. Assume that $X\times\Lambda\subseteq \Gamma$ is a local flow and $Z$ is a locally compact space. Let $\Pi_X: X\times \Lambda\rightarrow X$ and $\Pi_\Lambda:X\times \Lambda\rightarrow \Lambda$ be the canonical projection maps. See \cite{S} and \cite{Fr3}.

\begin{definition}
A parametrization of a local flow $X\subseteq \Gamma $ is a homeomorphism $\phi :Z\times \Lambda \rightarrow X$ such that for each $\lambda \in  \Lambda$, $\phi (Z\times \{\lambda\})$ is a local flow.
\end{definition}

Let $\phi:Z\times \Lambda \rightarrow X$ be a parametrization of a local flow $X$. Denote the restriction $\phi |_{(Z\times \{\lambda\})}$ by $\phi_\lambda$ and its image by $X_\lambda$.

\begin{lemma}\emph{[Salamon]}
For any compact set $N\subseteq X$ the set $\Lambda(N)=\{\lambda\in\Lambda \ |\ N\times \lambda$ is an isolating neighborhood in $X\times \lambda\}$ is open in $\Lambda$.
\end{lemma}

\begin{definition}
The space of isolated invariant sets is
$$
\mathscr{S}= \mathscr{S}(\phi)=\{S\times \lambda\ |\ \lambda\in \Lambda\ \text{and}\ S\times \lambda\ \text{is an isolated invariant compact set in}\ X\times \lambda\}.
$$
\end{definition}

For all compact sets $N\subseteq X$ define the maps $\varrho_N:\Lambda (N)\rightarrow \mathscr{S}$ and $\varrho_N(\lambda)=Inv(N\times \lambda)$. Then consider the topology on the space $\mathscr{S}$ generated by the sets $\{ \varrho_N(U)\ |\ N\subseteq X $ {compact}, $U\subseteq \Lambda(N)$ {open} $\}$.

A map $\gamma:\Lambda\rightarrow \mathscr{S}$ is called a \textit{section} of the space of isolated invariant sets if $\Pi_\Lambda\circ\gamma=id|_\Lambda$.

{We are interested in the situation where the homology index braids of admissible orderings of Morse decompositions at parameters $\lambda$ and $\mu$ are isomorphic, see Theorem \ref{teo_trancas}. That is, it is not enough that a Morse decomposition continues over $\Lambda$, it must also continue with a partial order, more specifically:}

\begin{figure}[!h]
  \centering
  \includegraphics[width=10cm]{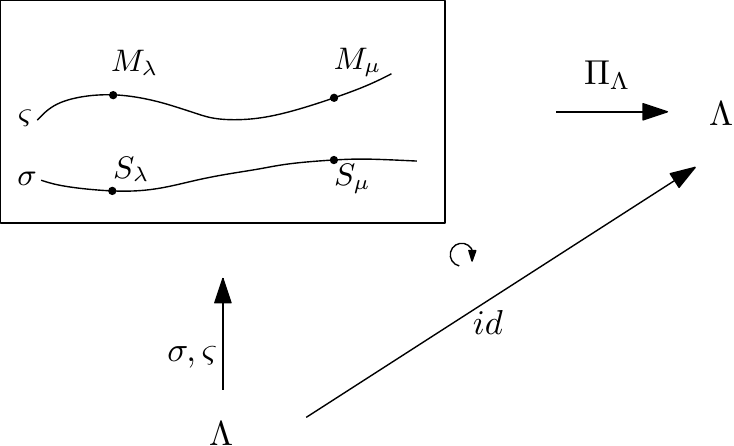}\\
  \caption{Sections from Definition \ref{def_continuation}}\label{secao}
\end{figure}

\begin{definition}\label{def_continuation}
Let $\mathcal {M}(S)=\{M(\pi)\ |\ \pi\in(\PP,<)\}$ be an ordered Morse decomposition of the isolated invariant set $S\subseteq X\times \Lambda$. Let $M_\lambda=\{M_\mu(\pi)\}_{\pi\in (\PP,<_\lambda)}$, $M_\mu=\{M_\mu(\pi)\}_{\pi\in (\PP,<_\mu)}$, $S_\lambda$ and $S_\mu$ be the sets obtained by intersection of $\mathcal{M}(S)$ and $S$ by the fibers $X\times\{\lambda\}$ and $X\times\{\mu\}$, respectively, where $<_\nu$ is the order restricted to the order $<$ in the parameter $\nu\in\Lambda$.

\begin{itemize}
  \item We say that $\mathcal{M}(S)$ with its order $<$ continues over $\Lambda$ if there exist sections $\sigma$ and $\varsigma_\pi : \Lambda \rightarrow \mathscr{S}$ such that $\{ \varsigma_\pi(\nu)|\ \pi \in (\PP,<_\nu)\}$ is a Morse decomposition for $\sigma(\nu),\ \forall \nu\in\Lambda$.
  \item  If, furthermore, there exist a path $\omega:[0,1]\rightarrow \Lambda$ from $\lambda$ to $\mu$; $\sigma(\lambda)=S_{\lambda}$; $\sigma(\mu)=S_{\mu}$; $\varsigma_\pi(\lambda)=M_{\lambda}(\pi)$; $\varsigma_\pi(\mu)=M_{\mu}(\pi)$; and if $\mathcal{M}(S)$ continues at least over $\omega([0,1])$, then we say that the admissible orderings $<_\lambda$ and $<_\mu$ are related by continuation or continue from one to the other. See Figure \ref{secao}.
\end{itemize}
\end{definition}

The following Lemma \ref{lemma_juncao} is a combination of Lemma 6.4 in \cite{S} and of Proposition 2.9 in \cite{MM1}.

\begin{lemma}\label{lemma_juncao}\emph{[McCord, Mischaikov, Salamon]}
\begin{itemize}
  \item Let $\gamma:\Lambda\rightarrow\mathscr{S}$ be a section, then $\gamma$ is continuous if and only if
  $$
  S=\displaystyle\bigcup_{\lambda\in\Lambda}\gamma(\lambda)
  $$
  is an isolated invariant set in $X\times \Lambda$.

  \item Let
  $$
  S=\displaystyle\bigcup_{\lambda\in\Lambda}\sigma(\lambda),\ \ \ M(\pi)=\displaystyle\bigcup_{\lambda\in\Lambda}\varsigma_\pi(\lambda)\ \text{for any}\ \pi\in \PP.
  $$
  Then, $S$ is an isolated invariant set in $X\times \Lambda$ under $\phi$ and $\mathcal{M}(S)=\{M(\pi)\ |\ \pi\in (\PP,<)\}$ is its Morse decomposition if, and only if, $\mathcal{M}(S)$ with its order continues.
\end{itemize}
\end{lemma}

Note that, by Lemma \ref{lemma_juncao}, Definition \ref{def_continuation} is equivalent to the definitions of continuation with order presented in \cite{S}, \cite{Fr3}, \cite{MM1} and \cite{FM}.

By Lemma \ref{lemma_juncao}, the \textit{minimal order} $<_m$ for $\mathcal{M}(S)$ that continues over $\Lambda$ is the flow defined order for $\mathcal{M}(S)$. Note that if $\mathcal{M}(S)$ with order $<$ continues then $<$ extends $<_m$.

\begin{proposition}
If $p<_m q$ then there exists $s_1,\ s_2,$ $\ldots,$ $s_{n}\in [0,1]$ and a sequence $(p_i)\subseteq\PP$ such that $p_0=q,\ p_{n}=p$ and the set of connecting orbits $C\left(M_{\omega(s_i)}(p_{i-1}),M_{\omega(s_i)}(p_i)\right)$ is non-empty, where $\omega:[0,1]\rightarrow \Lambda$ is a path between $\lambda$ and $\mu$. We call these connections, unordered chain connections, in short, $ucc$.
\end{proposition}
\prooff
Since a Morse set of $\mathcal{M}(S)$ is
$$M(\pi)=\displaystyle\bigcup_{\lambda\in\Lambda}\varsigma_\pi(\lambda)$$
then $p<_m q$ implies that there is a sequence $(p_i)\subseteq\PP$ such that $p_0=q,\ p_{n}=p$ and the set of connecting orbits $C\left(M(p_{i-1}),M(p_i)\right)$ is nonempty. Note that the connecting orbit between $M(p_{i-1})$ and $M(p_i)$ occurs at some parameter in $\Lambda$. Therefore, we have the desired result whenever $\mathcal{M}(S)$ continues over a path $\omega:[0,1]\rightarrow \Lambda$.
\cqd

Now we have the necessary framework to state the following results in \cite{Fr3} and \cite{FM}, which we use subsequently.

\begin{theorem}\label{teo_trancas}\emph{[Franzosa]} If the admissible orderings $<_\lambda$ and $<_\mu$ are related by continuation, then $\mathcal{H}(<_\lambda)$ and $\mathcal{H}(<_\mu)$, the homology index braids of the admissible orderings, are isomorphic.
\end{theorem}

\begin{theorem}\emph{[Franzosa]} If the admissible orderings $<_\lambda$ and $<_\mu$ are related by continuation, then $\mathcal{CM}(<_\lambda)=\mathcal{CM}(<_\mu)$.
\end{theorem}

\begin{proposition}\emph{[Franzosa]}
Let $<_1$ and $<_2$ be admissible orderings for $\mathcal{M}(S)$ and assume that $<_1$ is an extension of $<_2$. Then $$\mathcal{CM}(<_2)\subseteq \mathcal{CM}(<_1).$$
\end{proposition}

\begin{theorem}\emph{[Franzosa]}\label{cont_flow_order}
There exists a neighborhood $W$ of $\lambda$ in $\Lambda$ such that if $\mu\in W$, then $M_\lambda$ is related by continuation with order to a Morse decomposition $M_\mu$ of an isolated invariant set in $X_\mu$, and for such $M_\mu$, $\mathcal{CM}(M_\mu)\subseteq \mathcal{CM}(M_\lambda)$.
\end{theorem}

\begin{proposition}\label{Tchain}\emph{[Franzosa, Mischaikow]}
Let $C = \{C(p)\}_{p\in\PP}$ and $C' = \{C'(p)\}_{p\in\PP}$ be collections of
graded modules, and $\Delta:C(\PP)\rightarrow C(\PP)$ and $\Delta':C'(\PP)\rightarrow C'(\PP)$ be $<$-upper triangular boundary maps. If $T : C(\PP) \rightarrow C'(\PP)$ is $<$-upper triangular and such that $T\Delta=\Delta'T$, then $\{T(\I)\}_{\I\in \I(<)}$ is a chain map from $\mathcal{C}\Delta$ to $\mathcal{C}\Delta'$.
\end{proposition}

\section{Generalization of the Topological Transition Matrix.}

Suppose that $S_0$ and $S_1$ are invariant sets related by continuation in $X_{\lambda_0}$ and $X_{\lambda_1}$. Hence, there exists a map $\omega:[0,1]\rightarrow \Lambda$ such that $\omega(0)=\lambda_0$ and $\omega(1)=\lambda_1$ and an isolated invariant set $S$ over $\omega(I)$ such that $S_{\lambda_i}=S_i$. The inclusion $f_i:X_{\lambda_i}\rightarrow X\times \omega(I)$ induces an isomorphism $CH_\ast(S_i)\stackrel{f_{i\ast}}{\longrightarrow} CH_\ast(S)$, where $CH_\ast(S_i)$ and $CH_\ast(S)$ indicates the Conley homology indices of $S_i$ in $X_{\lambda_i}$ and of $S$ in $X\times \omega(I)$, respectively. Thus, there is an isomorphism
$$
F_\omega:CH_\ast(S_0) \stackrel{f^{-1}_{1\ast}\circ f_{0\ast}}{\longrightarrow} CH_\ast(S_1)
$$
that depends on the endpoint-preserving homotopy class $\omega$. If $\pi_1(\Lambda)=0$ then $F_\omega$ is independent of the path $\omega$ and one writes $F_{\lambda_1,\lambda_2}$ instead of $F_\omega$. The flow-defined continuation isomorphism is well-behaved with respect to composition of paths: $F_{\lambda,\lambda}=id$, $F_{\mu,\nu}\circ F_{\lambda,\mu}=F_{\lambda,\nu}$ and $F_{\lambda,\mu}=F^{-1}_{\mu,\lambda}$. For more details see \cite{MM2} and \cite{S}.

Let $M_{\lambda}=\{M_{\lambda}(\pi)\}_{\pi\in\mathbf{P}}$ and $M_{\mu}=\{M_{\mu}(\pi)\}_{\pi\in\mathbf{P}}$ be Morse decompositions, related by continuation, for the isolated invariant sets $S_\lambda\subseteq X_\lambda$ and $S_\mu\subseteq X_\mu$, respectively. In this setting $\lambda$ and $\mu \in \Lambda'\subseteq \Lambda$ and
$$
\Lambda' = \left\{ \lambda\in\Lambda\ | \ S_\lambda=\bigcup_{p\in \PP} M_\lambda(p)\right\}
$$
is a parameter set in which the corresponding Morse decomposition does not have connecting orbits. {To simplify notation, we denote $CH_\ast(M_\nu(\I))=H_{\ast,\nu}(\I)$ or just $CH(M_\nu(\I))=H_{\nu}(\I)$, where $\I\in\I(<_\nu)$ and $\nu\in\{\lambda,\mu\}$}.

By Conley's theory we have that there is an isomorphism $\Phi_\lambda:C_\ast\Delta_\lambda(\PP) \rightarrow H_{\ast,\lambda}(\PP)$ for $\lambda\in\Lambda'$, where $C_\ast\Delta_\lambda(\PP)=\bigoplus_{\pi\in\PP}CH(M_\lambda(\pi))$ is the chain complex with connection matrix $\Delta_\lambda$.

Therefore, we can carry out the continuation along the path $\omega$ in two ways: first by continuing $S_\lambda$ along the path $\omega$ using the isomorphism $F_{\lambda,\mu}$; secondly continuing $\bigcup_{p\in \PP} M_\lambda(p)$ along the path $\omega$ by using isomorphism $E_{\lambda,\mu}=\bigoplus_{p\in \PP} F_{\lambda,\mu}(M(p))$. More precisely, we have the following diagram

$$
\xymatrixcolsep{3pc}\xymatrix{
C\Delta_\lambda(\PP) \ar[r]^{E_{\lambda,\mu}} \ar[d]^{\Phi_\lambda} &  C\Delta_\mu(\PP) \ar[d]^{\Phi_\mu}\\
H_{\lambda}(\PP) \ar[r]^{F_{\lambda,\mu}}  &  H_{\mu}(\PP)
}
$$

In general the diagram above is not commutative. Due to the lack of commutativity, one is able to obtain information about connection orbits. Fix a base $\mathfrak{B}_\lambda$ in $C_\ast\Delta_{\lambda}$ and use the isomorphism $E_{\lambda, \mu}$ in order to define a base $E_{\lambda, \mu}(\mathfrak{B}_\lambda)$ in $C_\ast\Delta_{\mu}$. The composition $T_{\lambda,\mu}=\Phi_\mu^{-1}\circ F_{\lambda,\mu}\circ \Phi_\lambda$ can be represented as a matrix with respect to those bases and such matrix is called a \textbf{\textit{topological transition matrix}}.

The following theorem in \cite{MM1} summarizes some important properties for this matrix, which we refer to, from now on, as the classical topological transition matrix.

\begin{theorem}\label{teo_classical_top}\emph{[McCord-Mischaikow]} Let $\Lambda'\subseteq \Lambda$ be such that for all $\lambda$ and $\mu\in\Lambda'$ there are no connection orbits in $M_\lambda$ and $M_\mu$, and $M_\lambda$ and $M_\mu$ are related by continuation. Then
\begin{description}
  \item[(i)] $\Delta_\mu T_{\lambda,\mu} + T_{\lambda,\mu}\Delta_\lambda = 0$;
  \item[(ii)] $T_{\lambda,\mu}$ is an isomorphism;
  \item[(iii)] $T_{\lambda,\mu}$ is upper triangular matrix with respect to order $<$ ;
  \item[(iv)] If $\nu\in\Lambda'$ then $T_{\lambda,\lambda}=id$, $T_{\lambda,\nu}=T_{\mu,\nu}\circ T_{\lambda,\mu}$ and $T_{\mu,\lambda}=T_{\lambda,\mu}^{-1}$;
  \item [(v)] If $T_{\lambda,\mu}(p,q)\neq0$ and $\omega$ is a path between $\lambda$ and $\mu$, then there exists a finite sequence $0<s_1\leq s_2\leq \ldots \leq s_{n}<1$ and a sequence $(p_i)\subseteq\PP$ such that $p_0=q,\ p_{n}=p$ and the connecting orbit set $C\left(M_{\omega(s_i)}(p_{i-1}),M_{\omega(s_i)}(p_i)\right)$ is nonempty.
\end{description}
\end{theorem}

Item (i) of the Theorem above is trivial since $\Delta_\lambda=\Delta_\mu=0$, and $T_{\lambda,\mu}$ has the property of being unique, since $T_{\lambda,\mu}$ is a composition of isomorphisms.

%The following proposition from [\ref{Fr3}] ensures that the diagram of the next definition is well defined for any adjacent pairs $(\I,\J)$.

%\begin{proposition}\emph{[{Franzosa}]}
%Sejam $M_\lambda=\{M_\mu(\pi)\}_{\pi\in \PP}$ e $M_\mu=\{M_\mu(\pi)\}_{\pi\in \PP}$ decomposi\c c\~oes de Morse ordenadas, e assuma que as ordens admiss\'iveis associadas s\~ao relacionadas por continua\c c\~ao. Ent\~ao os conjuntos de Morse $M_\lambda(\I)$ e $M_\mu(\I)$ s\~ao relacionados por continua\c c\~ao para cada $\I\in\I(<)$, e o pares atrator-repulsor $( M_\lambda(\I), M_\lambda(\J) )$ em $M_\lambda(\I\J)$ e $( M_\mu(\I), M_\mu(\J) )$ em $M_\mu(\I\J)$ s\~ao relacionados por continua\c c\~ao para cada $(\I,\J)\in \I_2(<)$.
%\end{proposition}

\begin{definition}\label{def_cover}
Given chain complex braids $\mathcal{C}$ and $\mathcal{C}'$ and graded module braids $\mathcal{G}$ and $\mathcal{G}'$, a chain map $\mathcal{T}:\mathcal{C}\rightarrow \mathcal{C}'$ is said to \textbf{cover} an isomorphism $\theta$ (relative to $\Phi$ and $\Phi'$) if for all $\I\in \mathcal{I}(<)$, we have that the following diagram commutes
$$
\xymatrixcolsep{3pc} \xymatrix{
\mathcal{HC}(\I) \ar[r]^{\mathcal{T}_\ast(\I)} \ar[d]_{\Phi(\I)} & \mathcal{HC}'(\I) \ar[d]^{\Phi'(\I)}\\
\mathcal{G}(\I) \ar[r]^{\theta(\I)}                   & \mathcal{G}'(\I)
}
$$
where $\mathcal{T}_\ast(\I)$ is the homology map induced by the chain map $\mathcal{T}(\I)$, and $\Phi:\mathcal{HC}\rightarrow \mathcal{G}$ and $\Phi':\mathcal{HC}'\rightarrow \mathcal{G}'$ are isomorphisms from the homology of the chain complex braid to the graded module braid.
\end{definition}

\begin{definition}\label{GTM}
If, in {Definition} \emph{\ref{def_cover}}, $\mathcal{C}$ and $\mathcal{C}'$ arise from connection matrices $\Delta:\bigoplus_\PP C(p)\rightarrow \bigoplus_\PP C(p),\ \Delta':\bigoplus_\PP C'(p)\rightarrow \bigoplus_\PP C'(p)$, respectively, and $\mathcal{T}$ arises from a matrix $T:\bigoplus_\PP C(p)\rightarrow \bigoplus_\PP C'(p)$ then $T$ is called a \textbf{generalized transition matrix} for $\Delta$ and $\Delta'$.
\end{definition}

In a forthcoming paper we will further explore properties of the generalized transition matrix and demonstrate how it generalizes all three transition matrices, namely, singular \cite{R}, topological \cite{MM1} and algebraic \cite{FM}. In this paper we prove generalizations of the definition and properties of the classical topological transition matrices. With this in mind we restrict Definition \ref{GTM} in order to obtain a new {and broader definition for a topological transition matrix. The advantage over the classical setting is that the requirement of no connections between the Morse sets at the end parameter values in a continuation is dropped.}

\begin{definition}\label{GTTM}
If $T$ is a generalized transition matrix that covers the flow-defined continuation isomorphism $F$, then we refer to $T$ as a \textbf{generalized topological transition matrix}.
\end{definition}
We have the following characterization result.
\begin{proposition}\label{def_new}
$T$ is a \textbf{generalized topological transition matrix} related to the connection matrices $(\Delta_\lambda,\ \Phi_\lambda)$ and $(\Delta_\mu,\ \Phi_\mu)$, if and only if
$$T:\displaystyle \bigoplus_{p\in \PP}CH_\ast (M_\lambda(p)) \rightarrow \displaystyle \bigoplus_{ p\in  \PP}CH_\ast (M_\mu( p))$$
is a zero degree map such that

\begin{itemize}
  \item $\{T(\I)\}_{\I\in\I(<)}$ is a chain map from $\mathcal{C}\Delta_\lambda$ to $\mathcal{C}\Delta_\mu$;
  \item the following diagram commutes
\end{itemize}
\begin{table}[h!]\centering
\begin{tikzpicture}
  \matrix (m) [matrix of math nodes, row sep=3em,
    column sep=0.5em]{
	& H\Delta_\lambda(\emph{\I}) & & H\Delta_\lambda(\emph{\I\J}) & & H\Delta_\lambda(\emph{\J}) & & H\Delta_\lambda(\emph{\I})\\
	 H_\lambda(\emph{\I})       &  & H_\lambda(\emph{\I\J})       & & H_\lambda(\emph{\J})       & & H_\lambda(\emph{\I})     & \\
	& H\Delta_\mu(\emph{\I}) & & H\Delta_\mu(\emph{\I\J}) & & H\Delta_\mu(\emph{\J}) & & H\Delta_\mu(\emph{\I}) \\
	 H_\mu(\emph{\I})       & & H_\mu(\emph{\I\J})       & & H_\mu(\emph{\J})       & & H_\mu(\emph{\I})      &\\};	
  \path[-stealth]
    (m-1-2) edge (m-1-4) edge  (m-2-1)
            edge [densely dotted] (m-3-2)
    (m-2-1) edge [-,line width=6pt,draw=white] (m-2-3)
            edge (m-2-3) edge node [left] {{\small $F(\emph{\I})$}} (m-4-1)
    (m-3-2) edge [densely dotted] (m-3-4)
            edge [densely dotted] (m-4-1)
    (m-4-1) edge (m-4-3)
    (m-1-4) edge [densely dotted] (m-3-4) edge (m-2-3) edge (m-1-6)
    (m-3-4) edge [densely dotted] (m-4-3) edge [densely dotted] (m-3-6)
    (m-2-3) edge [-,line width=6pt,draw=white] (m-2-5) edge [-,line width=3pt,draw=white] (m-4-3)
            edge (m-4-3) edge (m-2-5)
	(m-4-3) edge (m-4-5)
	(m-1-6) edge (m-2-5) edge [densely dotted] (m-3-6) edge node [above] {{\small $\Delta_\lambda (\emph{\J},\emph{\I})$}} (m-1-8)
	(m-2-5) edge [-,line width=30pt,draw=white] (m-2-7) edge [-,line width=3pt,draw=white] (m-4-5) edge (m-4-5) edge node [above] {{\small $\delta_\lambda (\emph{\J},\emph{\I})$}} (m-2-7)
	(m-3-6) edge [densely dotted] (m-4-5) edge [densely dotted] (m-3-8)
	(m-4-5) edge (m-4-7)
	(m-1-8) edge (m-2-7) edge node [right] {{\small $\hat T(\emph{\I})$}} (m-3-8)
	(m-2-7) edge [-,line width=3pt,draw=white] (m-4-7) edge (m-4-7)
	(m-3-8) edge (m-4-7);
\end{tikzpicture}
\caption{}\label{dia_prin}
\end{table}
for all adjacent pairs \emph{$(\I,\J)$}, where $\hat T(\cdot)$ is the induced homology map of $T(\cdot)$.
\end{proposition}

\prooff
%Para a induzida $\hat T$ da definição anterior estar bem definida, precisa de que $T\circ %\Delta_\lambda=\Delta_\mu \circ T$.
By Definition \ref{GTTM}, the diagram commuting in the transversal sections implies that the whole diagram commutes. This follows easily since the top and bottom diagrams commute by the definition of the connection matrix; the diagram in the background commutes because $T$ is a chain map, and lastly the diagram in the foreground commutes by the continuation of the homology index braid (Theorem \ref{teo_trancas}) which, in part, asserts that for all adjacent pairs $(\I,\J)$ the following diagram commutes
$$
\xymatrixcolsep{3pc}\xymatrix{
\cdots \ar[r] & H_\lambda (\textbf{I}) \ar[r] \ar[d]^{F(\I)} & H_\lambda (\textbf{IJ}) \ar[r] \ar[d]^{F(\I\J)} & H_\lambda (\textbf{J}) \ar[r]^{\delta_\lambda(\textbf{J},\I)} \ar[d]^{F(\J)}& H _\lambda(\textbf{I})\ar[r] \ar[d]^{F(\I)}& \cdots\\
\cdots \ar[r] & H_\mu (\I) \ar[r] & H_\mu (\I\J) \ar[r] & H_\mu (\J) \ar[r]^{\delta_\mu(\J,\I)} & H_\mu (\I)\ar[r]& \cdots
}
$$
\cqd

Denote GTTM$(<)$ as the set of all generalized topological transition matrices with the partial order $<$.

When there are no connections in the $\lambda$ and $\mu$ parameters, then $\Delta_\lambda =0 =\Delta_\mu$. Thus, the induced homology map $\hat T=T$ and by choosing an adjacent pair $(\I,\J)$ such that $\PP=\I\J$, then from Diagram \ref{dia_prin} we have that the following diagram commutes

$$
\xymatrix{
 H_{\lambda}(\PP) \ar[r] \ar[d]^{F(\PP)} &  C\Delta_\lambda(\PP) \ar[d]^{T(\PP)}\\
 H_{\mu}(\PP) \ar[r] &  C\Delta_\mu(\PP)
}
$$

Therefore, $T(\PP)$ is a classical topological transition matrix. Thus, generalized topological transition matrices encompass the classical topological transition matrices.

Although the next result is straightforward from Definition \ref{def_new}, it is worthwhile to emphasize its importance, since given two connection matrices related by continuation there exists some entries that are the same for the matrices $\Delta_\lambda$ and $\Delta_\mu$. More accurately,

\begin{proposition}\label{pivo1gen}
Let $p,q\in \PP$ such that either $p$ and $q$ are not related by order $<$ or $p<q$. If the pair $(\{p\},\{q\})$ is an adjacent pair, i.e., there is no $p'\in\PP$ such that $p<p'<q$, then $\Delta_{qp,\lambda}=T^{-1}(\{q\})\circ \Delta_{qp,\mu}\circ T(\{p\})$.
\end{proposition}
\prooff
By hypothesis $(\{p\},\{q\})$ is an adjacent pair, then by the definition of the generalized topological transition matrix we have that the following diagram commutes

$$
\xymatrixcolsep{5pc}\xymatrix{
H\Delta_\lambda(\{p\}) \ar[r]^{\Delta_\lambda(\{q\},\{p\})} \ar[d]^{\hat T(\{p\})}       & H\Delta_\lambda(\{q\}) \ar[d]^{\hat T(\{q\})}  \\
H\Delta_\mu(\{p\}) \ar[r]^{\Delta_\mu(\{q\},\{p\})}                               & H\Delta_\mu(\{q\})
}
$$
which can rewritten as
$$
\xymatrixcolsep{3pc}\xymatrix{
CH(M_\lambda(p)) \ar[r]^{\Delta_{qp,\lambda}} \ar[d]^{ T(\{p\})}       & CH(M_\lambda(q)) \ar[d]^{ T(\{q\})}  \\
CH(M_\mu( p)) \ar[r]^{\Delta_{qp,\mu}}                                 & CH(M_\mu( q))
}
$$
Since $T(\{p\})$ and $T(\{q\})$ are isomorphisms it follows that $\Delta_{qp,\lambda}=T^{-1}(\{q\})\circ \Delta_{qp,\mu}\circ T(\{p\})$.
\cqd

In this paper we do not address the existence of the generalized topological transition matrix, {even though Corollary \ref{delta0} and Theorem \ref{teo_MS} establish the existence in particular cases.} 

A primary goal in this section is to establish properties of the generalized topological transition matrix - including connecting orbit existence results - corresponding to those of the classical topological transition matrix. In the last two sections we present applications of the generalized topological transition matrix.

%Since the existence problem is related to the unification of the transition matrix theory we will explore the existence question in a forthcoming paper.

The following properties of generalized topological transition matrices is an extension of Theorem \ref{teo_classical_top}.

\begin{theorem}\label{teo_prop}
Let $M_\lambda=\{M_\lambda(\pi)\}_{\pi\in (\PP,<_\lambda)}$ and $M_\mu=\{M_\mu(\pi)\}_{\pi\in (\PP,<_\mu)}$ be Morse decompositions, $\Delta_\lambda$ and $\Delta_\mu$ their respective connection matrices with the flow-defined order. Moreover, assume that $M_\lambda$ and $M_\mu$ are related by continuation with an admissible ordering $<$. Then the generalized topological transition matrix $T$ satisfies the following properties:
\begin{description}
  \item[(i)] $T \circ \Delta_\lambda=\Delta_\mu\circ  T$;
  \item[(ii)] $T_{\lambda,\mu}(\{p\})=id$ and $T$ is upper triangular with respect to $<$;
  \item[(iii)] $T$ is an isomorphism;
  \item[(iv)] $T_{\lambda,\lambda}=id$, $T_{\lambda,\nu}(\I)=T_{\mu,\nu}\circ T_{\lambda,\mu}(\I)$ and $T_{\mu,\lambda}(\I)=T_{\lambda,\mu}^{-1}(\I)$ are generalized topological transition matrices, for all intervals $\I\in \cal{I}$ and $p\in\PP$, in particular $T=T(\PP)$.
  \item [(v)] Let $\omega:[0,1]\rightarrow \Lambda$ be a path that continues $M_\lambda$ and $M_\mu$. Assume that $T_{\lambda,\mu}(p,q)\neq0$ for all generalized topological transition matrices and GTTM$(<_m)\neq \emptyset$ for all $\omega [s,t]$, where $s, t\in [0,1]$. Then there exists a finite sequence $0\leq s_1\leq s_2\leq \ldots \leq s_{n}\leq 1$ and a sequence $(p_i)\subseteq\PP$ such that $p_0=q,\ p_{n}=p$ and the set of connecting orbits $C\left(M_{\omega(s_i)}(p_{i-1}),M_{\omega(s_i)}(p_i)\right)$ is non-empty.

%  \begin{itemize}
%    \item if $T_{\lambda,\mu}(p,q)\neq 0$ and $\Delta_\lambda(\K)= 0 = \Delta_\mu(\K)$, where $\K$ is the interval that contains $p$ and $q$ at the ends;
%    \item or if $T_{\lambda,\mu}(p,q)\neq0$ for all generalized topological connection matrix and $\Delta_\lambda(\K)\neq 0$ or $\Delta_\mu(\K)\neq 0$\textbf{(Still need to prove)};
%  \end{itemize}

\end{description}
\end{theorem}

We prove Theorem \ref{teo_prop} below, but first we present a corollary.

\begin{corollary}\label{delta0}
Under the same hypothesis as Theorem \emph{\ref{teo_prop}}, assume that $T_{\lambda,\mu}(p,q)\neq 0$ and $\Delta_\lambda(\K)= 0 = \Delta_\mu(\K)$, where $\K$ is an interval that contains $p$ and $q$. Then there exist a finite sequence $0\leq s_1\leq s_2\leq \ldots \leq s_{n}\leq 1$ and a sequence $(p_i)\subseteq\PP$ such that $p_0=q,\ p_{n}=p$ and the set of connecting orbits $C\left(M_{\omega(s_i)}(p_{i-1}),M_{\omega(s_i)}(p_i)\right)$ is non-empty.
\end{corollary}
\prooff

Since $\Delta_\lambda(\K)= 0 = \Delta_\mu(\K)$ then
$$H\Delta_\lambda(\K)= \bigoplus_{\pi\in\K} CH(M_\lambda(\pi))\ \ \text{and}\ \ H\Delta_\mu(\K)= \bigoplus_{\pi\in\K} CH(M_\mu(\pi)),$$ hence GTTM$(<_m)=\{T_{\lambda,\mu}(\K)=\Phi_\mu^{-1}\circ F \circ \Phi_\lambda (\K)\} \neq \emptyset$. Therefore, by item (v) of Theorem \ref{teo_prop}, one just needs that $T_{\lambda,\mu}(p,q)\neq 0$ in order to prove the result.
\cqd

Note that when $\Delta_\lambda (\K)=0=\Delta_\mu(\K)$ we can not use the classical topological matrix to obtain Corollary \ref{delta0}, since $\Delta_\lambda (\K)=0$ does not imply that there is no connection at the parameter $\lambda$. Actually one can use Corollary \ref{delta0} whenever there is no connection at parameter $\lambda$ and $\mu$ in order to prove item (\textbf{v}) of Theorem \ref{teo_classical_top}.
\\

\prooff Items (\textbf{i}), (\textbf{ii}), (\textbf{iii}) and (\textbf{iv}) of Theorem \ref{teo_prop}.

(\textbf{i}) Since $\Delta_\lambda(\I)$ and $\Delta_\mu(\I)$ are boundary maps and $T_{\lambda.\mu}(\I)$ is a chain map, we have that $T_{\lambda.\mu}(\I)\circ \Delta_\lambda(\I) = \Delta_\mu(\I) \circ T_{\lambda.\mu}(\I)$ for all $\I$.

(\textbf{ii}) Fix a base $\mathfrak{B}_\lambda$ for the domain then  $\mathfrak{B}_\mu=\displaystyle \bigoplus_{p\in \PP} F(p)(\mathfrak{B}_\lambda)$ is a base for the codomain of the map
$$\displaystyle \bigoplus_{p\in \PP} F(p):\displaystyle \bigoplus_{p\in \PP}H_\lambda (p) \rightarrow \displaystyle \bigoplus_{ p\in \PP}H_\mu (p).
$$
Therefore, $\Phi_\lambda^{-1}(\mathfrak{B}_\lambda)$ and $\Phi_\mu^{-1}(\mathfrak{B}_\mu)$ are bases for the domain and codomain of the map $\bigoplus_{p\in \PP}  T_\ast(\{p\})$, i.e., $ T_\ast(\{p\})=id$ and, since $T_\ast(\{p\})= T(\{p\})$, then $T(\{p\})=id$.

In order to prove that $T$ is upper triangular it is enough to prove $T_{q,p}=0$ for $q\nless p$. %Indeed, $q\nless p$ implies that $p$ and $q$ are not related or $p<q$. If $p$ and $q$ are not related then there is no connection between $M_\nu(p)$ and $M_\nu(q)$, $\forall \nu\in\Lambda$, since $<$ is an extension of the flow defined order on $\mathcal{M}(S)$. Therefore, by $T(\{p,q\})$ being the classical topological transition matrix for the interval $\{p,q\}$, we have that $T(\{p,q\})=id$.

%Only it remains to prove for $p<q$,

Indeed let $\I$ be an interval that has $p$ and $q$ at the ends, and choose the adjacent pair $(p,\I\backslash p)$. It follows that
  $$
  \begin{array}{ccl}
    [T (\I)(\alpha\oplus 0)]&=& T_\ast (\I)[\alpha\oplus 0] = T_\ast (\I)\circ i (\alpha)= i\circ  T_\ast (\{p\})(\alpha)=\\
    &=&  i\circ T(\{p\})(\alpha) = i(\dot \alpha)=[i(\dot \alpha)]=[\dot\alpha\oplus 0],
  \end{array}
  $$
  where $\alpha \in C(p)$, $\dot \alpha=T(\{p\})(\alpha)$ and $0\in \bigoplus_{\pi\in\I\backslash p}C(\pi)$.
  Therefore, $$T(\I)(\alpha\oplus 0)-\dot\alpha\oplus 0\in Im\Delta_\mu(\I)$$ and since $\Delta'(\I)$ is an upper triangular boundary map then there exists $\beta\in\bigoplus_\I C'(\pi)$ such as
  $$T(\alpha\oplus 0)-\dot\alpha\oplus 0 = \Delta'(\I)(\beta)=d \oplus 0,$$
  where $d \in \bigoplus_{\pi\in\I\backslash q}C'(\pi)$.

  As $T(\alpha\oplus 0)=\dot\alpha\oplus\cdots\oplus T_{q,p}\cdot\alpha$ then $T_{q,p}\cdot\alpha=0$, thus $T_{q,p}=0$.

(\textbf{iii}) By item (ii) we have that $T$ is {an upper triangular matrix with nonzero entries on the diagonal}, thus $T$ is an isomorphism.

(\textbf{iv}) By the properties of the continuation isomorphism, specifically $F_{\lambda,\lambda}=id$, $F_{\mu,\nu}\circ F_{\lambda,\mu}=F_{\lambda,\nu}$ and $F_{\lambda,\mu}=F^{-1}_{\mu,\lambda}$, the result follows.
\cqd

One needs a technical lemma in order to prove item (\textbf{v}) of Theorem \ref{teo_prop}.

\begin{lemma}\label{lemma_ucc} Let $M_\lambda=\{M_\lambda(\pi)\}_{\pi\in (\PP,<_\lambda)}$ and $M_\mu=\{M_\mu(\pi)\}_{\pi\in (\PP,<_\mu)}$ be Morse decompositions, $\Delta_\lambda$ and $\Delta_\mu$ their respective connection matrices with the flow-defined order. Moreover, assume that $M_\lambda$ and $M_\mu$ are related by continuation with an admissible ordering $<$.
If $T_{\lambda,\mu}(p,q)\neq0$ for all generalized topological transition matrices and GTTM$(<_m)\neq \emptyset$ then there is a $ucc$ from $q$ to $p$.
\end{lemma}

\prooff

Since $T(p,q)\neq 0$ for all $T\in $ GTTM($<$), then  $T(p,q) \neq 0$ for all $T\in$ GTTM($<_m$), given that $<_m$ is the minimal order that continues. Therefore, by $T$ being $<_m$-upper triangular, we have that $p<_m q$.\cqd

\prooff Item (\textbf{v}) of Theorem \ref{teo_prop}.

By Lemma \ref{lemma_ucc} we have $T_{pq}\neq 0$ for all $T\in$ GTTM($<$) implies that there exists a $ucc$ from $q$ to $p$.

\begin{figure}[!h]
  \centering
  \includegraphics[width=12cm]{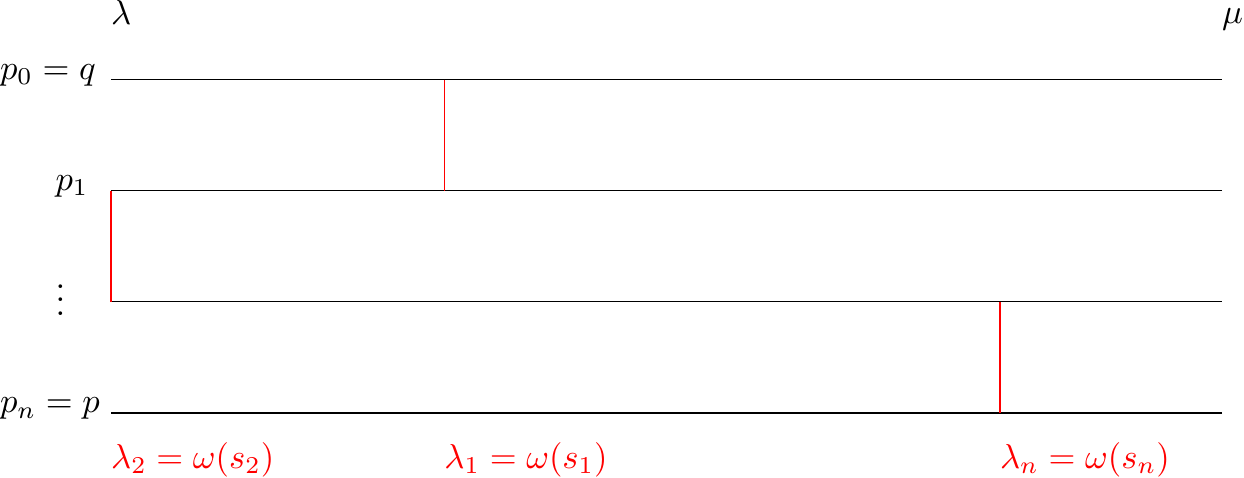}\\
  \caption{unordered chain connections from $q$ to $p$}\label{ccd1}
\end{figure}

It remains to prove that $s_i\leq s_{i+1}$, i.e., the connections occur along the path $\omega$. Since GTTM($<_m$)$\subseteq$GTTM($<$) then it is enough to prove the result for $T_{\lambda \mu}(p,q)\neq 0$ for the order $<_m$.

Indeed, we will prove by induction on the numbers of elements $m$ that $\K$ has.

Case $m=2$: Follows directly from Lemma \ref{lemma_ucc}.

Assume the result is true for $k<m+1$.

\begin{figure}[!h]
  \centering
  \includegraphics[width=15cm]{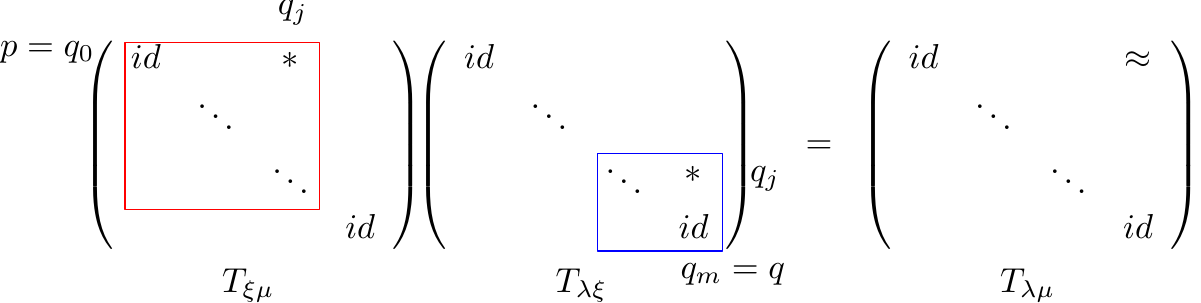}
  \caption{}\label{Comp_matriz3}
\end{figure}

Let $0<s_\xi<1$ and $\xi=\omega(s_\xi)$. Since $T_{\lambda \mu}(p,q)\neq 0$ for all $T_{\lambda \mu}\in\text{GTTM}_{\lambda \mu}(<)$, it follows that $\forall\ T_{\xi \mu}\in\text{GTTM}_{\xi\mu}(<)$ and $\forall\ T_{\lambda \xi}\in\text{GTTM}_{\lambda \xi}(<)$ there exists $0\leq j \leq m$ such that
$$
T_{\xi \mu}({q_0,q_j})\cdot T_{\lambda \xi}({q_j,q_m})\neq 0.
$$
See Figure \ref{Comp_matriz3}.

Since $j$ depends on the choice of $T_{\xi \mu}$ and $T_{\lambda \xi}$, we will fix $j$ for $T_{\xi \mu}\in\text{GTTM}_{\xi\mu}(<_{m_{\xi \mu}})$ and $T_{\lambda \xi}\in\text{GTTM}_{\lambda \xi}(<_{m_{\lambda \xi}})$, given that GTTM$_{\xi \mu}(<_{m_{\xi \mu}})\subseteq$GTTM$_{\xi \mu}(<_m)\subseteq$ GTTM$_{\xi \mu}(<$) and GTTM$_{\lambda \xi}(<_{m_{\lambda \xi}})\subseteq$GTTM$_{\lambda \xi}(<$), where $<_{m_{\xi \mu}}$ and $<_{m_{\lambda \xi}}$ are the minimal orders that continue for $\omega[s_\xi,1]$ and for $\omega[0,s_\xi]$, respectively.

\begin{figure}[!h]
  \centering
  \includegraphics[width=12cm]{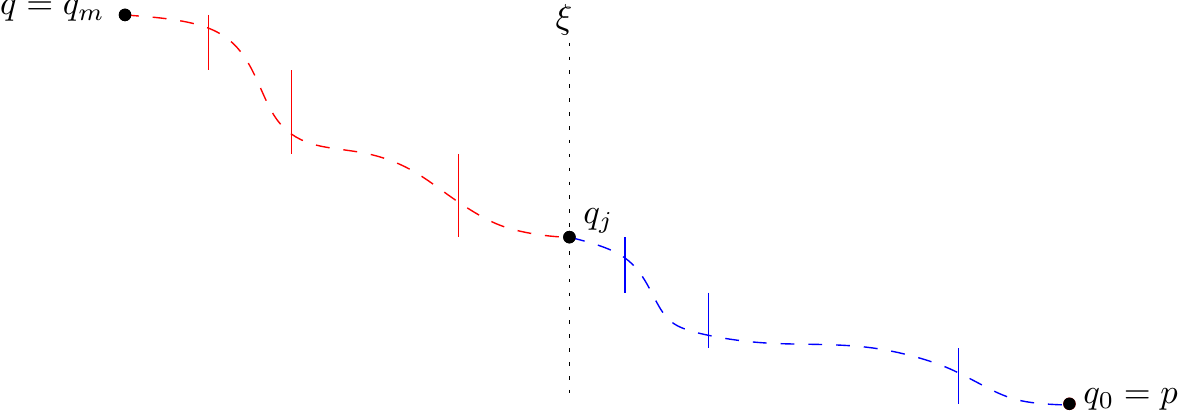}
  \caption{Connections between $M_\lambda(q_m)$ and $M_\mu(q_0)$ with $s$'s ordered.}\label{ccd7}
\end{figure}

If $j\neq 0$ and $j\neq m$ consider the submatrices $(T_{\xi \mu}({q_l,q_k}))_{l,k\in [0,j]}$ and $(T_{\lambda \xi}({q_l,q_k}))_{l,k\in [j,m]}$ in Figure \ref{Comp_matriz3}.
By the induction hypothesis there exist connections between $M_{\lambda}(q_m)$ and $M_{\xi}(q_j)$ and between $M_{\xi}(q_j)$ and $M_{\mu}(q_0)$ with $s$'s ordered. Using those connections we get connections between $M_\lambda(q_m)$ and $M_\mu(q_0)$ with $s$'s ordered. See Figure \ref{ccd7}.

\begin{figure}[!h]
  \centering
  \includegraphics[width=12cm]{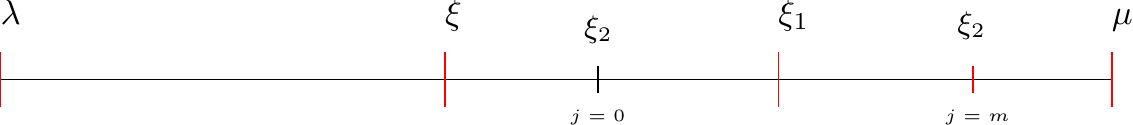}
  \caption{}\label{ccd8}
\end{figure}

If $j=0$ or $j=m$ then $T_{\lambda \xi}({p,q})\neq 0$ or $T_{\xi \mu}({p,q})\neq 0$, respectively. Suppose $T_{\xi \mu}({p,q})\neq 0$ then let $s_\xi < s_{\xi_1} < 1$ and $\xi_1=\omega(s_{\xi_1} )$. If for $\xi$ and $\mu$ there exists $j_{\xi_1}\neq 0$ and $j_{\xi_1}\neq m$ the result follows. Nevertheless if $j_{\xi_1}=0$ ($T_{\xi \xi_1}(p,q)\neq 0$) choose $s_\xi < s_{\xi_2} < s_{\xi_1}$ and $\xi_2=\omega(s_{\xi_2} )$ and if $j_{\xi_1}=m$ $(T_{\xi_1 \mu}(p,q)\neq 0)$ choose $s_{\xi_1}<s_{\xi_2}<1$ and $\xi_2=\omega(s_{\xi_2} )$. See Figure \ref{ccd8}.

By respecting this process and assuming always that $j_{\xi_i}=0$ or $j_{\xi_i}=m$, it follows that $\theta$ exists such that
$$\omega(s_{\xi_{i}})=\xi_i\rightarrow \theta,\ \omega(s_{\xi_{i'}})=\xi_{i'}\rightarrow \theta,\ 0\leq s_{\xi_i}<s_{\xi_{i'}}\leq 1 \text{ and } T_{\xi_i \xi_{i'}}(p,q)\neq 0.$$

%Changing this ``Thus $T_{\theta\theta}\neq id$ which is a contradiction. Therefore the result follows." to:

Suppose that there is no chain of connections from $q$ to $p$ at the parameter $\theta$, thus $p\nless_\theta q$, where $<_\theta$ is the flow defined order at parameter $\theta$. By Theorem \ref{cont_flow_order} there exists a neighborhood $W$ of $\theta$ such that $\mathcal{M}(S)$ continues with order $<_\theta$ over $W$, therefore $<_\theta$ extends $<_{m_W}$, where $<_{m_W}$ is the minimal order that continues over $W$. Choose $\xi_l$ and $\xi_{l'}$ such that $\theta\in \omega[s_{\xi_l},s_{\xi_{l'}}]\subset W$. It follows that $T_{\xi_l \xi_{l'}}(p,q)\neq 0$, i.e., there exists a $ucc$ from $q$ to $p$ in $\omega[s_{\xi_l},s_{\xi_{l'}}]$. Thus $p<_{m_W}q$. Since $<_\theta$ extends $<_{m_W}$ we have that $p<_\theta q$ for the Morse decomposition $\mathcal{M}(S)$, which is a contradiction, given that $p\nless_\theta q$. Therefore the result follows.
%The case $\Delta_\lambda(\K)\neq 0$ or $\Delta_\mu(\K)\neq 0$ is basically the same argument used in the previous case, with the following modifications:
%\begin{itemize}
%  \item The equalities and the inequalities must be true for all generalized topological transition matrix;
%  \item The compositions are true for the induced maps, i.e., it may happen that \\ $T_{\lambda_1 \mu}\circ T_{\lambda_2 \lambda_1} \circ T_{\lambda \lambda_2}\neq T_{\lambda \mu}$ but $\widehat{(T_{\lambda_1 \mu}\circ %T_{\lambda_2 \lambda_1} \circ T_{\lambda \lambda_2})}=\hat T_{\lambda \mu}$;
%  \item The marked entries, e.g. $\ast$-entries, will always be equal to $0$ or equal to $\approx$ for all generalized topological transition matrix.
%\end{itemize}
\cqd

The following Example \ref{example_Fr} shows that the connections obtained from item (v) of Theorem \ref{teo_prop} can occur at parameter $\lambda$ or $\mu$. This contrasts Theorem \ref{teo_classical_top} in \cite{MM1} and Theorem 3.13 in \cite{R}, where those connections can not occur at parameters $\lambda$ and $\mu$. However, when they do occur, one sees in Example \ref{example_Fr} that the classical results on transition matrices (topological, singular and algebraic) do not apply and therefore do not provide information on connections.

\begin{example}\label{example_Fr}Consider the following family of ordinary differential equations parameterized by the variable $\theta>0$:
$$
\dot x =y, \ \ \ \dot y=-\theta y - x(x-\frac{1}{3})(1-x).
$$

The connection matrices for $\theta>0$ are well known, see \cite{Fr2}, \cite{Fr3} and \cite{R}. Let $\mu$ be the parameter which has a heteroclinic connection between the Morse sets $M_\mu(2)$ and $M_\mu(3)$, and $0<\lambda<\mu$. The order that continues is the total order and the set of connection matrices are
$$
\left(\begin{array}{ccc}
0&\approx&0\\
0&0&0\\
0&0&0\\
\end{array}\right)
\ \ \ \text{and}\ \ \ \
\left(\begin{array}{ccc}
0&\approx&\approx\\
0&0&0\\
0&0&0\\
\end{array}\right).
$$

In this case, it is easy to see that $id\in \mathcal{T}^U_{\lambda,\ \mu}$, where $\mathcal{T}^U_{\lambda,\ \mu}$ is the set of algebraic transition matrices. Hence, Theorem 4.3, in \cite{FM}, does not apply and one can not retrieve dynamical information from $\mathcal{T}^U_{\lambda,\ \mu}$.

Now, in order to calculate the singular transition matrix, first introduce a slow drift $$\dot{\theta}=\epsilon (\lambda-\theta)(\mu-\theta)$$ in the parameter space, see Figure \ref{drift}.

\begin{figure}[!h]
  \centering
  \includegraphics[scale=0.7]{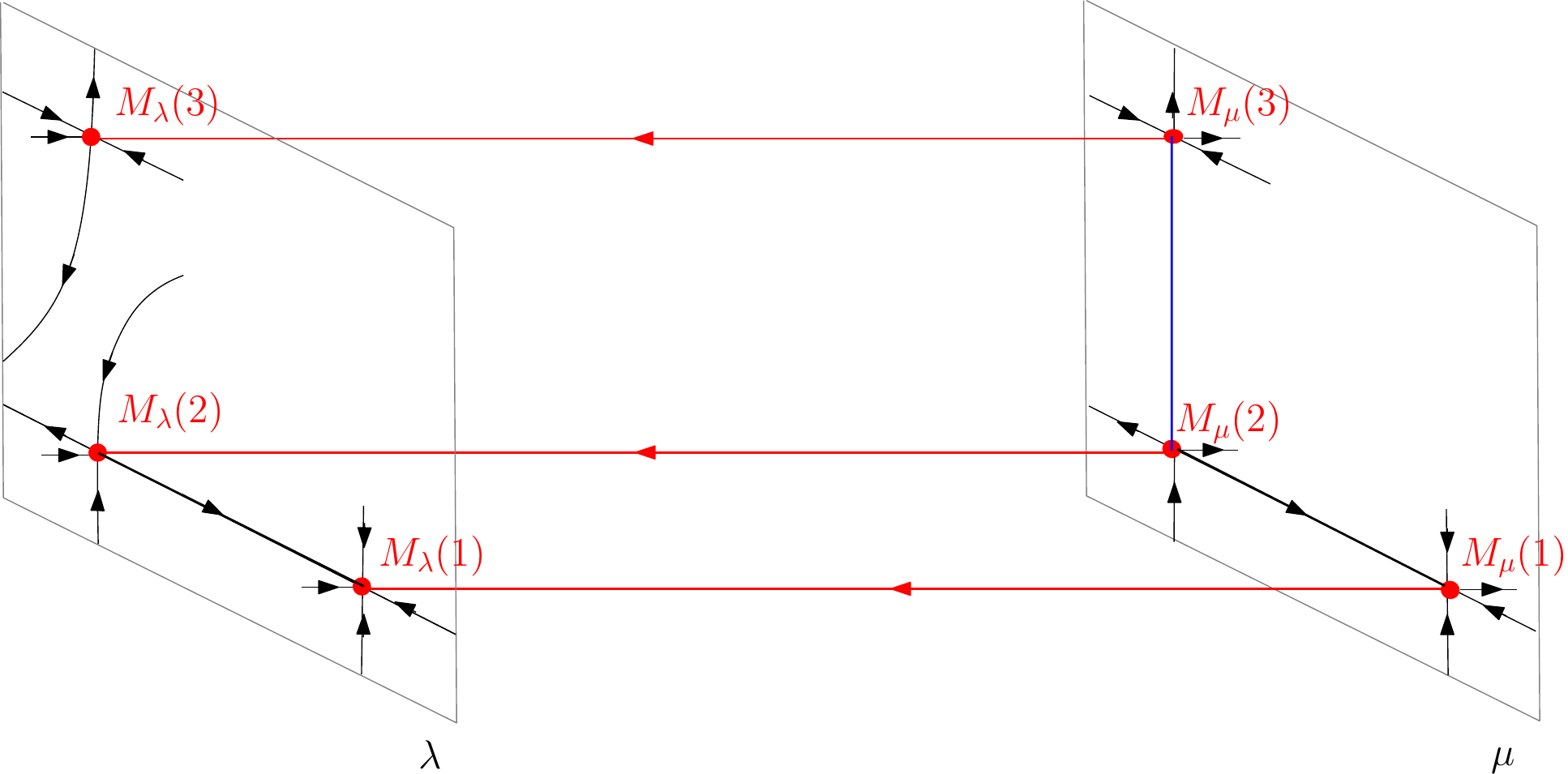}
  \caption{Drift flow in the parameter space.}\label{drift}
\end{figure}

Since there are two connection matrices for the flow ordering at parameter $\mu$, choose
$$
\Delta_\mu=
\left(\begin{array}{ccc}
0&\approx&\approx\\
0&0&0\\
0&0&0\\
\end{array}\right).
$$

The associated singular transition matrix is
$$
\Delta=\left(\begin{array}{cc}
\Delta_{\lambda} & T_s\\
0 & \Delta_{\mu}
\end{array}\right)
=\left(\begin{array}{cccccc}
0&\approx&0&\approx&&\\
0&0&0&&\approx&\ast\\
0&0&0&&&\approx\\
0&0&0&0&\approx&\approx\\
0&0&0&0&0&0\\
0&0&0&0&0&0\\
\end{array}\right).
$$

Since $\Delta^2=0$ must be zero, we see that multiplying the top row and the last column in $\Delta$ forces $\ast=\Delta(M_\mu(3),M_\lambda(2))\neq 0$. Even though the entry $T_s(p,q)\neq 0$, one can not conclude from $T_s$ that there exists a connection between $M_\lambda(2)$ and $M_\mu(3)$ for $\epsilon>0$ sufficiently small. Thus, one can not apply Theorem 3.13 in \cite{R}. Therefore, in this example, singular transition matrices do not give dynamical information. The same thing happens for classical topological transition matrices, since they are only defined when there are no connections at parameters $\lambda$ and $\mu$.

Now, we will calculate a generalized topological transition matrix related with
$$
\Delta_{\lambda}=
\left(\begin{array}{ccc}
0&\approx&0\\
0&0&0\\
0&0&0\\
\end{array}\right)
\ \ \ \text{and}\ \ \
\Delta_\mu=
\left(\begin{array}{ccc}
0&\approx&\approx\\
0&0&0\\
0&0&0\\
\end{array}\right).
$$
It is not hard to see that
$$
T_{\lambda,\mu}=
\left(\begin{array}{ccc}
\approx&0&0\\
0&\approx&\approx\\
0&0&\approx\\
\end{array}\right)
$$
satisfies item (i), (ii) and (iii) of Theorem \ref{teo_prop}, hence $T_{\lambda,\mu}$ is a good candidate to be a generalized topological transition matrix. Thus, it remains to prove that $T_{\lambda,\mu}$ covers the flow defined isomorphism $F$. Indeed, it is enough to prove that $T_{\lambda,\mu}(\I)$ covers $F(\I)$ for the interval $\I=\{2,3\}$ ($T_{\lambda,\mu}$ covers $F$ for others intervals, because $T_{\lambda,\mu}$ is a $0$ degree chain map). For $\I=\{2,3\}$ we have that $\Delta_\lambda(\I)=0=\Delta_\mu(\I)$, therefore we can choose generators: $\alpha_\nu$ and $\beta_\nu$ for $CH_1(M_\nu(\I))$; $a_\nu$ for $CH_1(M_\nu(3))$; and $b_\nu$ for $CH_1M_\nu(2)$, where $\nu=\lambda, \mu$. Figure \ref{CCI} indicates such choices.

\begin{figure}[!h]
  \centering
  \includegraphics[scale=0.9]{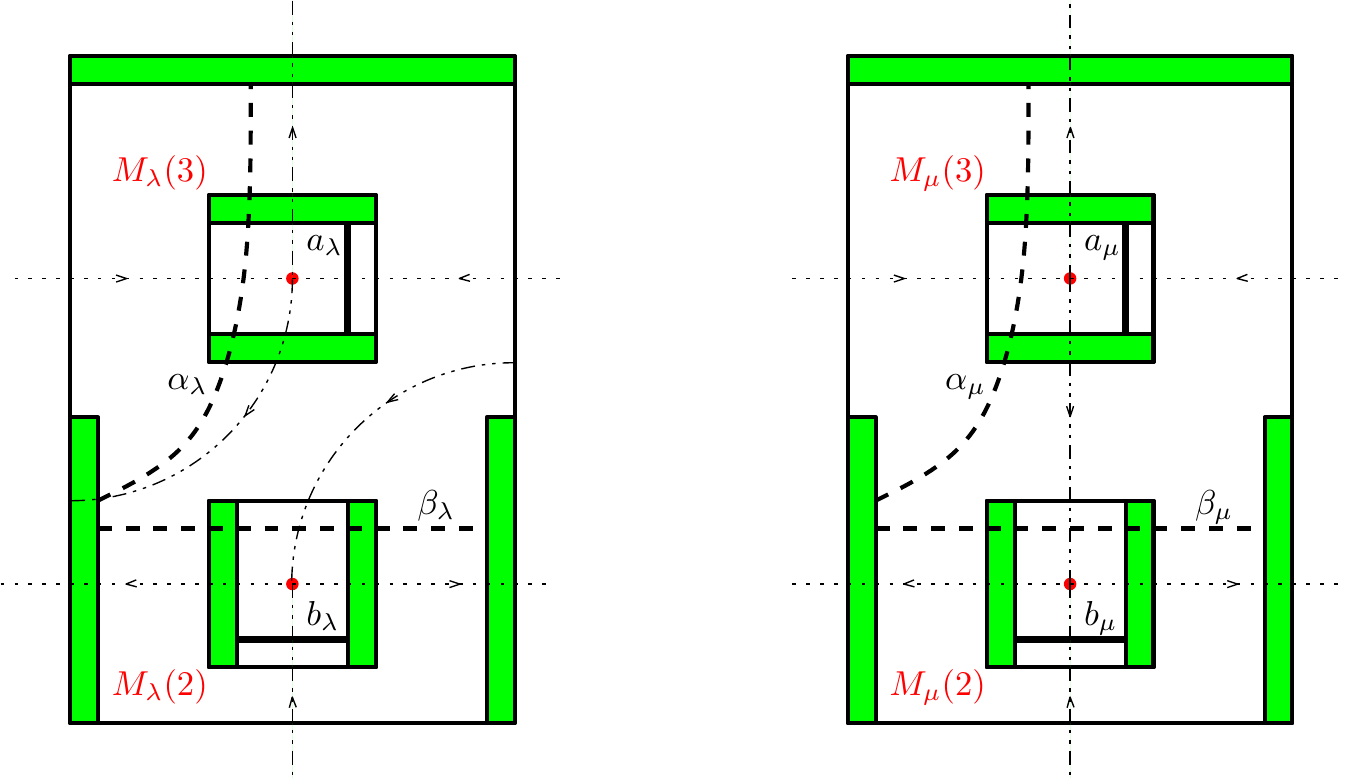}
  \caption{Generators of the Conley indices.}\label{CCI}
\end{figure}

Observe that $F_{\lambda, \mu}(\I)(\alpha_\lambda)=\alpha_\mu$, $F_{\lambda, \mu}(\I)(\beta_\lambda)=\beta_\mu$, $\Phi_\lambda(a_\lambda)=\alpha_\lambda$, $\Phi_\lambda(b_\lambda)=\beta_\lambda$, $\Phi_\mu(a_\mu)=\alpha_\mu \ast \beta_\mu$ and $\Phi_\mu(b_\mu)=\beta_\mu$. Thus $T_{\lambda,\mu}(\I)=\Phi_\mu^{-1}\circ F_{\lambda,\mu}\circ \Phi_\lambda(\I)$, which means that $T_{\lambda,\mu}(\I)$ covers $F(\I)$ and $T_{\lambda,\mu}(2,3)\neq 0$.

Note that the total order $1<2<3$ is the minimal order that continues, therefore by Lemma \ref{lemma_ucc} we have a connection between $M(2)$ and $M(3)$, since $T_{\lambda,\mu}(2,3)\neq 0$. In fact, one could have used item (v) of Theorem \ref{teo_prop}, since it is not hard to obtain that GTTM$(<_m)\neq \emptyset$ for all $[s,t]$, where $s, t\in [\lambda,\mu]$.

\end{example}

%No item (\textbf{v}) do theorem anterior, no caso em que $\K=(p,q)$ e $\Delta_\lambda(\K)\neq 0$ ou $\Delta_\lambda(\dot \K)\neq 0$ não precisa analisar $T_{\lambda,\mu}(p,q)$, porque a seguinte proposição garante a conexão
%
%\begin{proposition}
%Se $\Delta_{pq}\neq 0$ e $\I=(p,q)$ então $C(M_\lambda(p),M_\lambda(q))\neq \emptyset$.
%\end{proposition}
%
%\dem
%Note que $(\{p\},\{q\})$ é par adjacente, então segue, pela matriz de conexão, que
%$$
%\begin{array}{ccccc}
%\cdots \rightarrow &           CH(M_\lambda(q)) & \stackrel{\Delta_\lambda(\{p\},\{q\})}{\longrightarrow} &           CH(M_\lambda(p)) & \rightarrow \cdots\\
%                   &    \downarrow      &                                                 &   \downarrow       &                   \\
%\cdots \rightarrow &           H_\mu( q) & \stackrel{\delta_\lambda(\{ p\},\{ q\})}{\longrightarrow} &           H_\mu( p) & \rightarrow \cdots .
%\end{array}
%$$
%Por hipótese $\Delta_{pq}\neq 0$, daí $\delta_{pq}\neq 0$. Logo pela proposição 4.3 de [\ref{Fr1}] temos $C(M_\lambda(p),M_\lambda(q))\neq \emptyset$. \cqd

\section{Morse-Smale flows without periodic orbits}

%LEMBRA de colocar a frase: Connection Matrix for Morse-Smale flows without periodic orbits is unique for the flow defined order, for GTTM it is unique also.

In this section, the generalized topological transition matrix for Morse-Smale flows without periodic orbits is presented. In other words, the Morse decomposition consists of hyperbolic rest points and whenever the stable manifold of $M(\pi)$ and the unstable manifold of $M(\pi')$ have nonempty intersection, it is transversal.

As one can see in \cite{R2}, the connection matrix for Morse-Smale flows without periodic orbits is unique for the flow-defined order. It is no surprise that the generalized topological transition matrix is unique. This is verified in Theorem \ref{teo_MS}.

Furthermore in \cite{Mc} and \cite{S2}, an alternative and easier way to compute the connection matrix in this setting is presented. Likewise, we show in Theorem \ref{teo_MS} that the generalized topological transition matrix can be computed, without difficulty, from the set of the classical topological transition matrix.

%{More specifically, the connection matrix $\Delta$ is the boundary operator of the Morse complex, where $\Delta(x,y)$ counts the number of connecting orbits ``with orientation" between  critical points of consecutive indices $x$ and $y$.}

\begin{theorem}\label{teo_MS}
Let $M_\lambda=\{M_\mu(\pi)\}_{\pi\in \PP}$ and $M_\mu=\{M_\mu(\pi)\}_{\pi\in \PP}$ be Morse decompositions, and $\Delta_\lambda$ and $\Delta_\mu$ be the respective connection matrices with the flow-defined order. Moreover, assume that $M_\lambda$ and $M_\mu$ are related by continuation with the admissible ordering $<$ and the flow at $\lambda$ and $\mu$ is Morse-Smale without periodic orbits. Then the generalized topological transition matrix $T$ satisfies the following properties:
\begin{description}
  \item[(i)] $T\circ \Delta_\lambda =  \Delta_\mu \circ T$;
  \item[(ii)] $T_{\lambda,\mu}(\{p\})=id$ and $T$ is upper triangular with respect to $<$;
  \item[(iii)] $T$ is an isomorphism;
  \item[(iv)] $T_{\lambda,\lambda}=id$, $T_{\lambda,\nu}(\I)=T_{\mu,\nu}\circ T_{\lambda,\mu}(\I)$ and $T_{\mu,\lambda}(\I)=T_{\lambda,\mu}^{-1}(\I)$, for all intervals $\I\in \cal{I}$.
  \item [(v)] Let $\omega:[0,1]\rightarrow \Lambda$ be the path that continues $M_\lambda$ and $M_\mu$. If
   $T_{\lambda,\mu}(p,q)\neq0$ then there exist a finite sequence $0< s_1\leq s_2\leq \ldots \leq s_{n}< 1$ and a sequence $(p_i)\subseteq\PP$ such that $p_0=q,\ p_{n}=p$ and the set of connecting orbits $C\left(M_{\omega(s_i)}(p_{i-1}),M_{\omega(s_i)}(p_i)\right)$ is non-empty.
  \item [(vi)] $T_{\lambda,\mu}$ is unique;
  \item [(vii)] The generalized topological transition matrix is a matrix in block form with submatrices being the classical topological transition matrix $T_{top,i}$ of the critical points of index $i$
  $$T_{\lambda,\mu}=
\left(
  \begin{array}{cccc}
    T_{top,0} & 0         & 0      & 0 \\
    0         & T_{top,1} & 0      & 0 \\
    0         & 0         & \ddots & 0 \\
    0         & 0         & 0      & T_{top,k} \\
  \end{array}
\right),
  $$
\end{description}
\end{theorem}

\prooff
(\textbf{vii})
We let $Crit_i$ denote the set of critical points of index $i$, and $Crit_{i,i+1}$ denote the set of critical points that have index $i$ or $i+1$ along with all connecting orbits between them.

Since the order $<$ extends the flow ordering, then, without loss of generality, we can suppose the columns of the generalized topological transition matrix $T$ are ordered from the critical point of the lowest index to the largest index.

If $M(\I) \subseteq Crit_k$ then $T(\I)=T_{top}(\I)$, since there is no connecting orbit between critical points with the same index at parameters $\lambda$ and $\mu$. Moreover, for $M_\lambda(p)$ and $M_\lambda(q)$ with different indices, we have that $T(p,q)=0$, since $T$ is a zero degree map. Therefore $T$ must be a matrix in block form, and from property (ii) of Theorem \ref{teo_prop}, $T$ must be an upper triangular matrix. Thus
$$T=
\left(
  \begin{array}{cccc}
    T_{top,0} & 0         & 0      & 0 \\
    0         & T_{top,1} & 0      & 0 \\
    0         & 0         & \ddots & 0 \\
    0         & 0         & 0      & T_{top,k} \\
  \end{array}
\right)
$$
is a good candidate to be a generalized topological transition matrix.

Therefore, we will prove that:
\begin{itemize}
\item $\{T(\I)\}_{\I\in\I(<)}$ is a chain map from $\mathcal{C}\Delta_\lambda$ to $\mathcal{C}\Delta_\mu$;
\item $T$ makes Diagram \ref{dia_prin} commute for all adjacent pairs $(\I,\J)$.
\end{itemize}
Indeed, by Proposition \ref{Tchain} we only need to show that $T\circ \Delta_\lambda =  \Delta_\mu \circ T$, i.e.,
$$
\left(
  \begin{array}{cccc}
    T_{top,0} & 0         & 0      & 0 \\
    0         & T_{top,1} & 0      & 0 \\
    0         & 0         & \ddots & 0 \\
    0         & 0         & 0      & T_{top,k} \\
  \end{array}
\right)
\left(
  \begin{array}{cccc}
    0 & \Delta_\lambda(Crit_1,Crit_0)        & 0      & 0 \\
    0         & 0 & \ddots      & 0 \\
    0         & 0         & 0 & \Delta_\lambda(Crit_k,Crit_{k-1}) \\
    0         & 0         & 0      & 0 \\
  \end{array}
\right)
=$$
$$=
\left(
  \begin{array}{cccc}
    0 & \Delta_\mu(Crit_1,Crit_0)        & 0      & 0 \\
    0         & 0 & \ddots      & 0 \\
    0         & 0         & 0 & \Delta_\mu(Crit_k,Crit_{k-1}) \\
    0         & 0         & 0      & 0 \\
  \end{array}
\right)
\left(
  \begin{array}{cccc}
    T_{top,0} & 0         & 0      & 0 \\
    0         & T_{top,1} & 0      & 0 \\
    0         & 0         & \ddots & 0 \\
    0         & 0         & 0      & T_{top,k} \\
  \end{array}
\right).
$$
By multiplying the matrices one obtains
$$
T_{top,l-1}\circ\Delta_\lambda(Crit_l,Crit_{l-1})=\Delta_\mu(Crit_l,Crit_{l-1})\circ T_{top,l}\ \ \text{for all}\ l\in\{1,...,k\}.
$$
But this follows from the commutativity of the following diagram
$$
\xymatrixcolsep{7pc}\xymatrix{
H\Delta_\lambda(Crit_l) \ar[r]^{\Delta_\lambda(Crit_l,Crit_{l-1})} \ar[d]^{T(Crit_l)}&  H\Delta_\lambda(Crit_{l-1}) \ar[d]^{T(Crit_{l-1})}\\
H\Delta_\mu(Crit_l)     \ar[r]^{\Delta_\mu(Crit_l,Crit_{l-1})}     & H\Delta_\mu(Crit_{l-1})
}
$$
whereas $T(Crit_l)=T_{top,l}=\Phi_\mu^{-1}\circ F(Crit_l)\circ \Phi_\lambda$.

Now we will prove that $T$ makes Diagram \ref{dia_prin} commute for all adjacent pairs $(\I,\J)$. Indeed, suppose without loss of generality, $M_\I\subseteq Crit_{k-1,k}$ and $M_\J\subseteq Crit_{k,k+1}$, and let $L_j=L\cap Crit_j$ where $L=\I$ or $L=\J$ and $j\in\{k-1,k,k+1\}.$

First, we calculate the homologies from the long exact sequences that come from Diagram \ref{dia_prin}
$$
\begin{array}{c}
\cdots\rightarrow H_{k+1}\Delta (\I)\rightarrow H_{k+1}\Delta (\I\J)\rightarrow H_{k+1}\Delta (\J)\rightarrow H_{k}\Delta (\I)\rightarrow H_{k}\Delta (\I\J)\rightarrow\\
\rightarrow H_{k}\Delta (\J)\rightarrow H_{k-1}\Delta (\I)\rightarrow H_{k-1}\Delta (\I\J)\rightarrow H_{k-1}\Delta (\J)\rightarrow \cdots
\end{array}
$$
By the definition of connection matrices, we have to make the following calculations:

$(C_\ast\Delta(\I\J),\Delta(\I\J)):$
$$
0 \rightarrow \displaystyle\bigoplus_{\pi\in\J_{k+1}}CH(M(\pi)) \stackrel{\Delta}{\rightarrow} \displaystyle\bigoplus_{\pi\in Crit_{k}}CH(M(\pi)) \stackrel{\Delta}{\rightarrow} \displaystyle\bigoplus_{\pi\in\I_{k-1}}CH(M(\pi))\rightarrow 0,
$$

$(C_\ast\Delta(\I),\Delta(\I)):$
$$
0 \rightarrow 0 {\rightarrow} \displaystyle\bigoplus_{\pi\in \I_{k}}CH(M(\pi)) \stackrel{\Delta}{\rightarrow} \displaystyle\bigoplus_{\pi\in\I_{k-1}}CH(M(\pi))\rightarrow 0,
$$

$(C_\ast\Delta(\J),\Delta(\J)):$
$$
0 \rightarrow \displaystyle\bigoplus_{\pi\in\J_{k+1}}CH(M(\pi)) \stackrel{\Delta}{\rightarrow} \displaystyle\bigoplus_{\pi\in \J_{k}}CH(M(\pi)) {\rightarrow} 0\rightarrow 0.
$$
Therefore the homologies are
$$
H_k\Delta(\I)=Ker\Delta_k(\I),\ H_{k-1}\Delta(\I)=\displaystyle\frac{\displaystyle\bigoplus_{\pi\in \I_{k-1}} CH(M(\pi))}{Im\Delta_k(\I)}\ \text{and}\  H_j\Delta(\I)=0\ \forall j\neq k,k-1
$$
$$
H_{k+1}\Delta(\J)=Ker\Delta_{k+1}(\J),\ H_{k}\Delta(\J)=\displaystyle\frac{\displaystyle\bigoplus_{\pi\in \J_{k}} CH(M(\pi))}{Im\Delta_{k+1}(\J)}\ \text{and}\  H_j\Delta(\J)=0\ \forall j\neq k+1,k
$$
$$
H_{k+1}\Delta(\I\J)=Ker\Delta_{k+1}(\I\J),\ H_{k}\Delta(\I\J)=\displaystyle\frac{Ker \Delta_k}{Im\Delta_{k+1}}\ \text{and}\  H_{k-1}\Delta(\I\J)=\displaystyle\frac{\bigoplus_{\I_{k-1}} CH(M(\pi))}{Im\Delta_{k}},
$$
for all $j\neq k-1,k,k+1.$
Thus, making the substitutions one obtains

\begin{table}[!h]
  \centering
  \includegraphics[width=16cm]{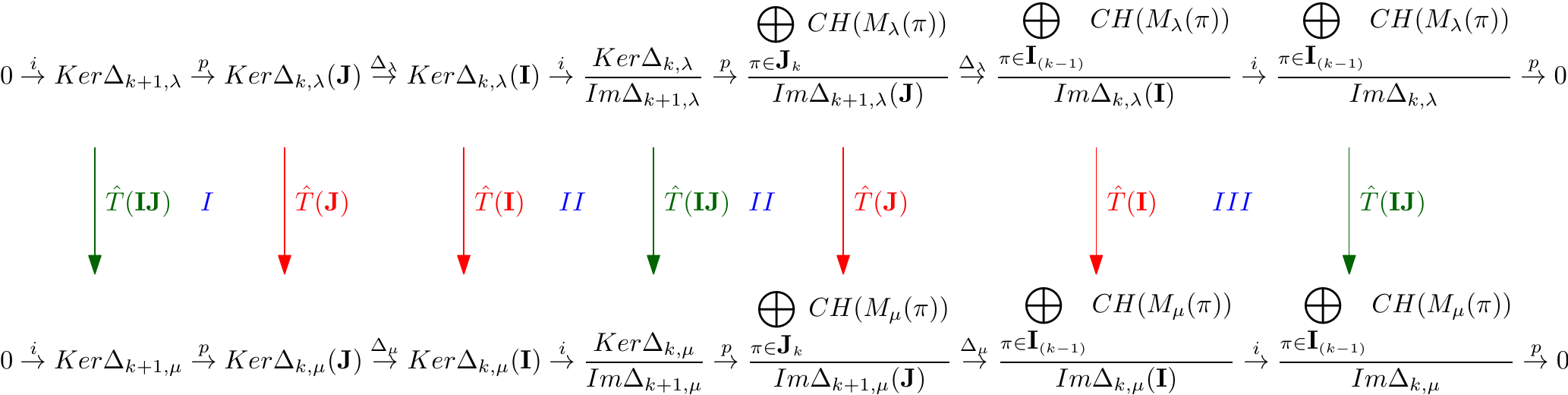}
 \caption{}\label{Seq1}
\end{table}

Now we will prove that $\textcolor[rgb]{1.00,0.00,0.00}{\hat T(\I)}$ and $\textcolor[rgb]{1.00,0.00,0.00}{\hat T(\J)}$ are the induced maps for the map $T$. Indeed, consider the adjacent pair $(\I_{k-1},\I_k)$, hence
$$
\xymatrixcolsep{2pc}\xymatrix{
  0 \ar[r]^<<<<{i} & H_k\Delta_\lambda (\I)   \ar[r]^{p} \ar[d]^{\hat T(\I)} & H_k\Delta_\lambda(\I_k) \ar[r]^<<<<<{\Delta_\lambda} \ar[d]^{\hat T(\I_k)} & H_{k-1}\Delta_\lambda(\I_{k-1}) \ar[r]^{i} \ar[d]^{\hat T(\I_{k-1})}& H_{k-1}\Delta_\lambda(\I) \ar[r]^<<<<<{p} \ar[d]^{\hat T(\I)} & 0 \\
  0 \ar[r]^<<<<{i} & H_k\Delta_\mu (\I)   \ar[r]^{p} & H_k\Delta_\mu(\I_k) \ar[r]^<<<<<{\Delta_\mu} & H_{k-1}\Delta_\mu(\I_{k-1}) \ar[r]^{i} & H_{k-1}\Delta_\mu(\I) \ar[r]^<<<<<{p} & 0
}
$$

Note that $T(\I_k)$ and $T(\I_{k-1})$ are the classical topological transition submatrices of $T$, i.e., the {rightmost} diagram commutes in Diagram \ref{Diagrama_mini}.

\begin{table}[h!]\centering
\begin{tikzpicture}
  \matrix (m) [matrix of math nodes, row sep=1.5em,
    column sep=0em]{
    & H_k\Delta_\lambda(\I) &  & H_k\Delta_\lambda(\I_k) \\
    H_{k,\lambda}(\I) & & H_{k,\lambda}(\I_k) & \\
    & H_k\Delta_\mu(\I) & & H_k\Delta_\mu(\I_k) \\
    H_{k,\mu}(\I) & & H_{k,\mu}(\I_k) & \\};
  \path[-stealth]
    (m-1-2) edge (m-1-4) edge  (m-2-1)
            edge [densely dotted] (m-3-2)
    (m-1-4) edge node [right] {{\small $T(\I_k)$}} (m-3-4) edge (m-2-3)
    (m-2-1) edge [-,line width=6pt,draw=white] (m-2-3)
            edge (m-2-3) edge node [left] {{\small $F(\I)$}} (m-4-1)
    (m-3-2) edge [densely dotted] (m-3-4)
            edge [densely dotted] (m-4-1)
    (m-4-1) edge (m-4-3)
    (m-3-4) edge (m-4-3)
    (m-2-3) edge [-,line width=3pt,draw=white] (m-4-3)
            edge (m-4-3);
\end{tikzpicture}
\caption{}\label{Diagrama_mini}
\end{table}

In order to show the {leftmost} diagram in Diagram \ref{Diagrama_mini} commutes when we place $\hat T(\I)$, it is enough to show the background diagram commutes, whereas the rightmost and foreground diagrams commute.

Indeed, as $H_k\Delta_\lambda(\I)=Ker\Delta_{k,\lambda}(\I) \unlhd H_k\Delta_\lambda(\I_k)=\bigoplus_{\I_k}CH(M_\lambda(\pi))$ it follows that
$$
T(\I_k)\circ p=p\circ \hat T(\I),
$$
since $p$, the induced projection map $C_k\Delta_\lambda(\I)\rightarrow C_k\Delta_\lambda(\I_k),$ is the inclusion map.
Analogously, one can show
$$
\hat T(\I)\circ i=i\circ T(\I_{k-1}),
$$
since $$H_{k-1}\Delta_\lambda(\I)=\displaystyle\frac{\displaystyle\bigoplus_{\pi\in\I_{k-1}}CH(M_\lambda(\pi))}{Im\Delta_\lambda(\I)}=\displaystyle\frac{H_{k-1}\Delta_\lambda(\I_{k-1})}{Im\Delta_\lambda(\I)}.$$
The latter assertion follows since $i$, the induced inclusion map $C_{k-1}\Delta_\lambda(\I_{k-1})\rightarrow C_{k-1}\Delta_\lambda(\I),$ is the projection map.

In the same way, the above construction can be done for $\J=\J_k\J_{k+1}$ and the adjacent pair $(\J_k,\J_{k+1})$. Hence $\textcolor[rgb]{1.00,0.00,0.00}{\hat T(\I)}$ and $\textcolor[rgb]{1.00,0.00,0.00}{\hat T(\J)}$ are the map induced by $T$ and make Diagram \ref{dia_prin} commute sectionwise.

Now it remains to prove the same for $\textcolor[rgb]{0.00,0.59,0.00}{\hat T(\I\J)}$, using the same idea described previously. So it is enough to prove that the diagrams $\textcolor[rgb]{0.00,0.07,1.00}{I},\ \textcolor[rgb]{0.00,0.07,1.00}{II}$ and $\textcolor[rgb]{0.00,0.07,1.00}{III}$ in Diagram \ref{Seq1} commute. But this comes from $\Delta_\lambda(\I)$ and $\Delta_\lambda(\J)$ being submatrices of $\Delta_\lambda(\I\J)$. The only diagram that deserves special attention is diagram $\textcolor[rgb]{0.00,0.07,1.00}{II}$.

We will prove that diagram $\textcolor[rgb]{0.00,0.07,1.00}{II}$ commutes.

Indeed, let $a\in Ker\Delta_{k,\lambda}(\I)$. Then $ i(a)=a+Im\Delta_{k+1,\lambda}(\I\J)=[a+b]$, where $b$ is such that there exists $c\in \displaystyle\bigoplus_{\pi\in\J_{k+1}}CH(M_\lambda(\pi)$ such that $\Delta_{k+1,\lambda}(\I\J)(c)=b.$ Applying $\hat T(\I\J)$ we have
$$
\hat T(\I\J)([a+b])=[T(\I\J)(a+b)].
$$
On the other hand,
$$
i\circ \hat T(\I)(a)=\hat T(\I)(a)+Im\Delta_{k+1,\mu}(\I\J)=[T(\I)(a)+d],
$$
where $d\in Im\Delta_{k+1,\mu}(\I\J)$. Therefore, we need to show $[T(\I\J)(a+b)]=[T(\I)(a)+d]$.

Since $T(\I\J)a=T(\I)a$, for $a\in Ker\Delta_{k,\lambda}(\I)\subseteq \displaystyle\bigoplus_{\pi\in\I_k} CH(M_\lambda(\pi))$ it is enough to show
$$
T(\I\J)a + T(\I\J)b - T(\I)a - d= T(\I\J)b - d\in Im\Delta_{k,\mu}(\I\J).
$$

Indeed, $d\in Im\Delta_{k,\mu}(\I\J)$ so it is sufficient to prove that Diagram \ref{diaast} commutes.

\vspace{-0.1cm}
\begin{table}[h!]
$$\xymatrixcolsep{4pc}\xymatrix{
\displaystyle\bigoplus_{\pi\in\J_{k+1}}CH(M_\lambda(\pi)) \ar[r]^<<<<<<<<<{\Delta_{k+1,\lambda}{(\I\J)}} \ar[d]^{ T(\J_{k+1})} & \displaystyle\bigoplus_{\pi\in\I_k}CH(M_\lambda(\pi))\oplus\displaystyle\bigoplus_{\pi\in\J_k}CH(M_\lambda(\pi)) \ar[d]^{T(\I\J)}\\
 \displaystyle\bigoplus_{\pi\in \J_{k+1}}CH(M_\mu(\pi)) \ar[r]^<<<<<<<<<<{\Delta_{k+1,\mu}{(\I\J)}} & \displaystyle\bigoplus_{\pi\in\I_k}CH(M_\mu(\pi))\oplus\bigoplus_{\pi\in\J_k}CH(M_\mu(\pi))
}$$\vspace{-0.5cm}
\caption{}\label{diaast}
\end{table}

When considering the adjacent pair $(\I_k\J_k,\J_{k+1})$, we have that the diagram commutes in the sections, since $T(\J_{k+1})$ and $T(\I_k\J_k)$ are precisely the classical topological transition matrices.

Therefore the diagram
\begin{center}
\begin{tikzpicture}
  \matrix (m) [matrix of math nodes, row sep=1.5em,
    column sep=0em]{
    & \bigoplus_{\textbf{J}_{k+1}}CH(M_\lambda(\pi)) &  & \bigoplus_{\textbf{I}_k}CH(M_\lambda(\pi))\oplus\bigoplus_{\textbf{J}_k}CH(M_\lambda(\pi)) \\
    H_{k+1,\lambda}(\textbf{J}_{k+1}) & & H_{k,\lambda}(\textbf{I}_k\textbf{J}_k) & \\
    & \bigoplus_{ \textbf{J}_{k+1}}CH(M_\mu(\pi)) & & \bigoplus_{\textbf{I}_k}CH(M_\mu(\pi))\oplus\bigoplus_{\textbf{J}_k}CH(M_\mu(\pi)) \\
    H_{k+1,\mu}(\textbf{J}_{k+1}) & & H_{k,\mu}(\textbf{I}_k\textbf{J}_k) & \\};
  \path[-stealth]
    (m-1-2) edge node [above] {{\small $\Delta_{k+1,\lambda}({\textbf{IJ}})$}} (m-1-4) edge  (m-2-1)
            edge [densely dotted] (m-3-2)
    (m-1-4) edge node [right] {{\small $T(\textbf{I}_k\textbf{J}_k)$}} (m-3-4) edge (m-2-3)
    (m-2-1) edge [-,line width=6pt,draw=white] (m-2-3)
            edge (m-2-3) edge node [left] {{\small $F(\textbf{J}_{k+1})$}} (m-4-1)
    (m-3-2) edge [densely dotted] (m-3-4)
            edge [densely dotted] (m-4-1)
    (m-4-1) edge (m-4-3)
    (m-3-4) edge (m-4-3)
    (m-2-3) edge [-,line width=3pt,draw=white] (m-4-3)
            edge (m-4-3);
\end{tikzpicture}
\end{center}
in the background commutes. Observing that
$$
T_k(\I\J):\displaystyle\bigoplus_{\I_{k-1}}0\displaystyle\bigoplus_{\pi\in\I_k\J_k}CH_k(M_\lambda(\pi))\displaystyle\bigoplus_{\J_{k+1}}0 \displaystyle\stackrel{0\oplus T_k(\I_k\J_k)\oplus 0}{ \xrightarrow{\hspace*{2cm}}}
\displaystyle\bigoplus_{\I_{k-1}}0\displaystyle\bigoplus_{\pi\in\I_k\J_k}CH_k(M_\mu(\pi))\displaystyle\bigoplus_{\J_{k+1}}0
$$
thus Diagram \ref{diaast} commutes.

Hence it was proved that diagram $\textcolor[rgb]{0.00,0.07,1.00}{II}$ commutes on the left. One can prove analogously that it also commutes on the right.

Therefore $T$ makes Diagram \ref{dia_prin} commute for all adjacent pairs $(\I,\J)$.

Item (\textbf{vi}) follows from the fact that the classical topological transition matrices are unique, and hence $T$ is unique.

Items (\textbf{i}), (\textbf{ii}), (\textbf{iii}) and (\textbf{iv}) follow from $T$ being a generalized topological transition matrix and from the fact that $T$ is unique.

Item (\textbf{v}) follows from (\textbf{vii}) and from property (\textbf{v}) of Theorem \ref{teo_prop}.
\cqd

\section{Generalized Transition Matrix for the Sweeping Method.}

In this section we present an application of a generalized topological transition matrix in a continuation associated to a dynamical spectral sequence, see \cite{CdRS} and \cite{FdRS}. Our dynamical interpretation result implies the existence of connecting orbits in a fast-slow system ``going from $M(q)$ to $M(p)$" for a nontrivial entry on $T^r_{pq}$ associated to the spectral sequence.

Let $M$ be an $n$-dimensional compact Riemannian manifold, $f : M \rightarrow \mathbb{R}$ a Morse function that is Morse-Smale, and $\phi$ the gradient flow of $f$. Choose a finite Morse decomposition $\bigcup M(p)$, $p \in \PP =\{1,...,m\}$ such that there are
distinct critical values $c_p$ with $f^{-1}(c_p) \supset M(p)$. Then
$$
\{F_p\}_{p=1}^m=\{f^{-1}(-\infty, c_p+\epsilon)\}^m_{p=1}
$$
This defines an admissible ordering on $M$ called the \textit{filtration order}. In this case, each Morse set, $M(p)$, is a
non-degenerate singularity of the gradient flow $\phi$ and hence, {the Conley index of each Morse set is the homology of
a pointed $k$-sphere, where $k$ is the Morse index of the singularity $M(p)$}. We denote by $h_k^{(p)}$ the index $k$ singularity in $F_p\backslash F_{p-1}$.

In the case where each $M(p)$ is a non-degenerate singularity and the stable and unstable manifolds intersect
transversally, the connection matrix $\Delta$ associated to $\mathcal{D}(M)$ is unique (see [R1] and [R2]). It can also be defined
as the differential of the graded Morse chain complex $(C,\Delta)$, where $C$ is generated by the singularities and
graded by their indices, i.e., $C = Z_2<Crit f>$ and $\Delta$ is determined by the maps $\Delta_k:C_k\rightarrow C_{k-1}$ via
$$
\Delta_k(x)=\sum_{y\in Crit_{k-1}f}n(x,y)<y>.
$$
where $n(x,y)$ is the number of connecting orbits counted mod 2 for nondegenerate singularities $x$ and $y$ of indices
$k$ and $k -1$ respectively. We require that the columns of the matrix $\Delta$ are ordered such that $\Delta_k(F_pC_k)\subset F_{p-1}C_{k-1}$.

In [CdRS] and [MdRS] it is proved that a certain algorithm (called the sweeping method) applied to  a connection matrix $\Delta$ determines a spectral sequence $(E^r,d_r)$ of a filtered chain complex $(C,\Delta)$. To achieve this we apply the sweeping
method to the connection matrix. The sweeping method is an iterative process that, given a connection matrix,
generates a collection of connection matrices $\Delta^1,...,\Delta^F$ and transition matrices $T^1,..., T^F$. This method singles out important nonzero entries, namely primary pivots and change of basis pivots, of the $r$-th diagonal of $\Delta^r$, which are necessary to define a matrix $\Delta^{r+1}$.
At each step, $\Delta^{r+1}$ is obtained from $\Delta^{r}$ by a change of basis.

The change of basis matrices $T^r$, $r=1,...,F$, determined by the sweeping method algorithm are called \textit{transition matrices associated to the spectral sequence}.

\begin{proposition}\label{prop_mari}\emph{[FdRS]}
Each matrix $T^r$ associated to the spectral sequence satisfies the following properties:
\begin{enumerate}
\item $\Delta T+T\Delta'=0$;
\item $T$ is an isomorphism;
\item $T$ is an upper triangular matrix with respect to the filtration order.
\end{enumerate}
\end{proposition}

In \cite{FdRS} a dynamical interpretation of the spectral sequence $(E_r,d_r)$ is given in a setting of a fast-slow system flow,
$$
\begin{array}{lcl}
\dot{x} &=& f(x,y)\\
\dot{y} &=& \varepsilon (y-1)(y-2),
\end{array}
$$
where the sweeping method output of $n\times n$ connection matrices and transition matrices, $\Delta^1, T^1,....,T^{F-1},\Delta^F$,
reveals bifurcations that arise as a result of the nonzero entries of $T^r$.

In this article we consider a more general fast-slow system flow
$$
\begin{array}{lcl}
\dot{x} &=& f(x,y)\\
\dot{y} &=& \varepsilon g(x,y),
\end{array}
$$ in $M\times [0,F]$ with the following properties:

\begin{itemize}
    \item When $\varepsilon = 0 $ the parameterized system has an isolated invariant set $S_y$ for each $y$
    which continues over the interval $[0,F]$ (slow variable), and which has a Morse decomposition
    $\mathcal{D}(M)_y =
    \{M_\lambda(p)_y \ |\ p=1,...,n\}$ that also continues over the interval  $[0,F]$.
    Assume also that the order from the sweeping method continues.

    \item For each $r=0,1,...,F$, we have that $f(x,r)=f_r$ is a gradient function
    which comes from a Morse function whose stable and unstable manifolds intersect transversely.
    And $g(M_y(p),y)\neq 0$ for $y\in(0,F)$ and $p\in \PP$.

    \item For each $r=0,1,...,F$, the sweeping method connection matrices $\Delta^r$ are connection matrices of the Morse decomposition at parameter $r$.

    \item Lastly, the continuation of the Morse decomposition
    $\mathcal{D}(M)_y$ is such that Diagram \ref{mwbox} commutes.

\begin{table}[h!]\centering
\begin{tikzpicture}
  \matrix (m) [matrix of math nodes, row sep=3em,
    column sep=-1.5em]{
	& \displaystyle \bigoplus_{x\in Crit_{k+1}f_r} \mathbb{Z} x & & \displaystyle \bigoplus_{x\in Crit_{k}f_r} \mathbb{Z} x & & \displaystyle \bigoplus_{x\in Crit_{k-1}f_r} \mathbb{Z} x \\	
	 H(Crit_{k+1}f_r)       &  & H(Crit_{k}f_r)       & & H(Crit_{k-1}f_r)       & \\
	& \displaystyle \bigoplus_{x\in Crit_{k+1}f_{r+1}} \mathbb{Z} x & & \displaystyle \bigoplus_{x\in Crit_{k}f_{r+1}} \mathbb{Z} x & & \displaystyle \bigoplus_{x\in Crit_{k-1}f_{r+1}} \mathbb{Z} x  \\
	 H(Crit_{k+1}f_{r+1})       & & H(Crit_{k}f_{r+1})       & & H(Crit_{k-1}f_{r+1})       & \\};	
  \path[-stealth]
    (m-1-2) edge node [above] {{\small $\Delta_{k+1,k}^r$}}  (m-1-4) edge   (m-2-1)
            edge [densely dotted] node [left] {{\small }}  (m-3-2)
    (m-2-1) edge [-,line width=6pt,draw=white] (m-2-3)
            edge (m-2-3) edge node [left] {{\small $F_{k+1}^{r}$}} (m-4-1)
    (m-3-2) edge [densely dotted] (m-3-4)
            edge [densely dotted] (m-4-1)
    (m-4-1) edge node [above] {{\small $\delta_{k+1,k}^{r+1}$}} (m-4-3)
    (m-1-4) edge [densely dotted] (m-3-4) edge (m-2-3) edge node [above] {{\small $\Delta_{k,k-1}^r$}} (m-1-6)
    (m-3-4) edge [densely dotted] (m-4-3) edge [densely dotted] (m-3-6)
    (m-2-3) edge [-,line width=6pt,draw=white] (m-2-5) edge [-,line width=3pt,draw=white] (m-4-3)
            edge (m-4-3) edge (m-2-5)
	(m-4-3) edge node [above] {{\small $\delta_{k,k-1}^{r+1}$}} (m-4-5)
	(m-1-6) edge (m-2-5) edge node [right] {{\small $T_{k-1}^{r}$}} (m-3-6)
	(m-2-5)  edge [-,line width=3pt,draw=white] node [right] {{\small }} (m-4-5) edge (m-4-5)
	(m-3-6) edge (m-4-5);
%	        \node [below=3cm, align=flush center,text width=8cm] {Morse Box diagram.};
\end{tikzpicture}
\caption{Morse Box diagram.}\label{mwbox}
\end{table}

At each stage $r$ of the sweeping method, $F_k^r$'s are continuation isomorphisms,
$\Delta^r_{k+1,k}$ is a submatrix of the connection matrix
$\Delta^r$ at parameter $r$, and $T_k^r$ is a submatrix of $T^r$.

\end{itemize}

\begin{theorem}
Consider the fast-slow system defined previously. Then the sweeping method transition matrices are generalized topological transition matrices.
\end{theorem}

Therefore, the sweeping method transition matrices inherit all properties from Theorem \ref{teo_MS}.
  $$T_{0,1}=
\left(
  \begin{array}{cccc}
    T_{top,0} & 0         & 0      & 0 \\
    0         & T_{top,1} & 0      & 0 \\
    0         & 0         & \ddots & 0 \\
    0         & 0         & 0      & T_{top,k} \\
  \end{array}
\right)
  $$

\prooff

By Proposition \ref{Tchain}, the collection of the submatrices $\{T^r(\I)\}_{\I\in\I(<)}$ of the sweeping method transition matrices is a chain map from $\mathcal{C}\Delta_r$ to $\mathcal{C}\Delta_{r+1}$, since $T^r$ is an upper triangular matrix and $T^r\circ \Delta^r=\Delta^{r+1}\circ T^r$.

Now we will prove that $T^r$ is a matrix in block form as $T_{r,r+1}$ the generalized topological transition matrix for $\Delta^r$ and $\Delta^{r+1}$.

Indeed, since the boundary map $\Delta^r$ of a Morse complex is a connection matrix (see \cite{S2}) we have that for each adjacent pair of intervals ($\I,\J$)
the following diagram is commutative
$$
\xymatrixcolsep{3pc}\xymatrix{
\cdots \ar[r] & H\Delta^r (\textbf{I}) \ar[r] \ar[d]^{\Phi_r(\textbf{I})}  & H\Delta^r (\textbf{IJ})
\ar[r] \ar[d]^{\Phi_r(\textbf{IJ})} & H\Delta^r (\textbf{J}) \ar[r]^{\Delta^r
(\textbf{J},\I)} \ar[d]^{\Phi_r(\textbf{J})}  & H\Delta^r (\textbf{I})\ar[r] \ar[d]^{\Phi_r(\textbf{I})} & \cdots
\\
\cdots \ar[r] & H_r (\textbf{I}) \ar[r] & H_r (\textbf{IJ}) \ar[r] & H_r (\textbf{J}) \ar[r]^{\delta^r
(\textbf{J},\I)} & H_r (\textbf{I}) \ar[r] & \cdots
}
$$

Set $M_r(\textbf{I})=Crit_{k-1}{f_r}$ and
$M_r(\textbf{J})=Crit_{k}{f_r}$. One obtains $\Delta^r(\textbf{J},\I)=\Delta^r_{k,k-1}$. And since there are no connections in
$M_r(\textbf{I})$ and in $M_r(\textbf{J})$ it follows that
$$
H\Delta^r (\textbf{I})=
\left\{%
\begin{array}{ll}
    \bigoplus \Z x, & \hbox{n=k-1;} \\
    0, & \hbox{otherwise.} \\
\end{array}%
\right.
$$
$$
H\Delta^r (\textbf{J})=
\left\{%
\begin{array}{ll}
    \bigoplus \Z x, & \hbox{n=k;} \\
    0, & \hbox{otherwise.} \\
\end{array}%
\right.
$$
Therefore Diagram \ref{dia_prin}, in this case, is equal to Diagram \ref{aux}.

\begin{table}[h!]\centering
\begin{tikzpicture}
  \matrix (m) [matrix of math nodes, row sep=3em,
    column sep=-1.8em]{
	& \ \ \ \ 0 \ \ \ \  & & H_k\Delta^r(\textbf{IJ}) & & \bigoplus_{x\in \textbf{J} } \mathbb{Z}x & & \bigoplus_{x\in \textbf{I}} \mathbb{Z}x \\	
	 H_k(Crit_{k-1}f_r)       &  & H_k(Crit_{k-1,k}f_r)       & & H_k(Crit_{k}f_r)       & &  H_{k-1}(Crit_{k-1}f_r)     & \\	
	& \ \ \ \ 0 \ \ \ \ & & H_k\Delta^{r+1}(  \textbf{I}  \textbf{J}) & & \bigoplus_{x\in  \textbf{J}} \mathbb{Z}x & & \bigoplus_{x\in  \textbf{I} } \mathbb{Z}x \\
	 H_k(Crit_{k-1}f_{r+1})       & & H_k(Crit_{k-1,k}f_{r+1})       & & H_k(Crit_{k}f_{r+1})       & & H_{k-1}(Crit_{k-1}f_{r+1})      &\\};	
  \path[-stealth]
    (m-1-2) edge (m-1-4) edge  (m-2-1)
            edge [densely dotted] (m-3-2)
    (m-2-1) edge [-,line width=6pt,draw=white] (m-2-3)
            edge (m-2-3) edge node [left] {{\small $F^r_{k-1}({\I})$}} (m-4-1)
    (m-3-2) edge [densely dotted] (m-3-4)
            edge [densely dotted] (m-4-1)
    (m-4-1) edge (m-4-3)
    (m-1-4) edge [densely dotted] (m-3-4) edge (m-2-3) edge (m-1-6)
    (m-3-4) edge [densely dotted] (m-4-3) edge [densely dotted] (m-3-6)
    (m-2-3) edge [-,line width=6pt,draw=white] (m-2-5) edge [-,line width=3pt,draw=white] (m-4-3)
            edge (m-4-3) edge (m-2-5)
	(m-4-3) edge (m-4-5)
	(m-1-6) edge (m-2-5) edge [densely dotted] (m-3-6) edge node [above] {{\small $\Delta^r_{k,k-1}$}} (m-1-8)
	(m-2-5) edge [-,line width=30pt,draw=white] (m-2-7) edge [-,line width=3pt,draw=white] (m-4-5) edge (m-4-5) edge  (m-2-7)
	(m-3-6) edge [densely dotted] (m-4-5) edge [densely dotted] (m-3-8)
	(m-4-5) edge  node [above] {{\small $\delta^{r+1}$}} (m-4-7)
	(m-1-8) edge (m-2-7) edge node [right] {{\small $ T^r_{k-1}({\I})$}} (m-3-8)
	(m-2-7) edge [-,line width=3pt,draw=white] (m-4-7) edge (m-4-7)
	(m-3-8) edge (m-4-7);
\end{tikzpicture}
           \caption{}\label{aux}
\end{table}

Note that for homology dimensions different from $k$ and $k-1$, the sequence becomes

$$
\xymatrixcolsep{3pc}\xymatrix{
0 \ar[r] \ar[d]^{\Phi_r(\textbf{I})} & H_n\Delta^r (\textbf{IJ}) \ar[r] \ar[d]^{\Phi_r(\textbf{I}\J)} & 0 \ar[r]^0 \ar[d]^{\Phi_r(\textbf{J})} & 0 \ar[d]^{\Phi_r(\textbf{I})}\\
H_n(Crit_{k-1}f_r)\ar[r] & H_n(Crit_{k-1,k}f_r) \ar[r] & H_n(Crit_{k}f_r)\ar[r]^<<<<<<<{\delta^r(\textbf{J},\I)} & H_{n-1}(Crit_{k-1}f_r)
}
$$

Since $  H_n\Delta^r(\textbf{IJ})=0$ and $\Phi_r$ is an isomorphism, it follows that $H_n(\cdot)=0$.

Therefore $T^r$ has the same block structure as $T_{r,r+1}$ and by hypothesis each submatrix $T^r_k$ of $T^r$ is actually $T_{top,k}$. Thus, by item (\textbf{vi}) and (\textbf{vii}) of Theorem \ref{teo_MS}, we have that $T^r=T_{r,r+1}$. \cqd

Applying Proposition \ref{pivo1gen} to a fast-slow system as defined previously, we obtain the following corollary.

\begin{corollary}\label{pivo1gen2}
The entries from $\Delta^0$ which are preserved, independent of continuation, are the primary pivots of the initial stages of the sweeping method.
\end{corollary}

Note that these entries are the first step to define a new algebraic method in attempting to generalize the sweeping method, which is defined only for Morse-Smale flows without periodic orbits.

\end{document}